\documentstyle[12pt]{article}
\newcommand{\numero}[1]{
\addtocounter{section}{1}
\begin{center}{\bf \thesection .\
#1\vspace{-.1in}}\end{center}
\setcounter{subsection}{0}
\setcounter{lemma}{0}\indent}

\newcommand{\subnumero}[1]{
\pagebreak[1]\begin{center}{\em #1}\nopagebreak\end{center}
}

\newcommand{\eop}{\hfill $\Box$\vspace{.1in}}

\newenvironment{inset}{
\begin{list}{}{\setlength{\rightmargin}{.3in}
\setlength{\leftmargin}{.3in}\setlength{\topsep}{.08in}
\setlength{\listparindent}{.13in}\setlength{\itemindent}{.13in}}
\small\item}{
\end{list}\normalsize}

\newtheorem{lemma}{Lemma}[section]
\newtheorem{theorem}[lemma]{Theorem}

\newtheorem{corollary}[lemma]{Corollary}

\newtheorem{conjecture}{Conjecture}
\newtheorem{proposition}[lemma]{Proposition}

\newcommand{\semidirect}{\mbox{\hspace*{.12
cm}$\times$\hspace*{-.38
cm}\raisebox{0cm}{$\triangleright$}\hspace*{.28cm}}}

\newcommand{\cc}{{\bf C}}

\newcommand{\qq}{{\bf Q}}

\newcommand{\zz}{{\bf Z}}

\newcommand{\Oo}{{\cal O}}
\newcommand{\Dd}{{\cal D}}

\newcommand{\Bb}{{\cal B}}
\newcommand{\Mm}{{\cal M}}

\newcommand{\Xx}{{\cal X}}

\begin{document}

\section*{Algebraic (geometric) $n$-stacks}

Carlos Simpson

\bigskip

In the introduction of Laumon-Moret-Bailly (\cite{LaumonMB} p. 2) they
refer to a
possible theory of algebraic $n$-stacks:
\begin{inset}
Signalons au passage que Grothendieck propose d'\'elargir \`a son tour le cadre
pr\'ec\'edent en rempla\c{c}ant les $1$-champs par des $n$-champs (grosso modo,
des faisceaux en $n$-cat\'egories sur $(Aff)$ ou sur un site arbitraire) et il
ne fait gu\`ere de doute qu'il existe une notion utile de $n$-champs
alg\'ebriques \ldots .
\end{inset}
The purpose of this paper is to propose such a theory. I  guess that the main
reason why Laumon and Moret-Bailly didn't want to get into this theory was for
fear of getting caught up in a horribly technical discussion of $n$-stacks of
groupoids over a general site. In this paper we simply {\em assume} that a
theory of $n$-stacks of groupoids exists.  This is not an unreasonable
assumption, first of all because there is a relatively good substitute---the
theory of simplicial presheaves or presheaves of spaces
(\cite{Brown} \cite{BrownGersten} \cite{Joyal} \cite{Jardine1}
\cite{kobe} \cite{flexible})---which should be equivalent, in an appropriate
sense, to any eventual theory of $n$-stacks; and second of all because it seems
likely that a real theory of $n$-stacks of $n$-groupoids could be developped in
the near future (\cite{BreenAsterisque}, \cite{Tamsamani}).

Once we decide to ignore the technical complications involved in theories of
$n$-stacks, it is a relatively straightforward matter to generalize Artin's
definition of algebraic $1$-stack. The main observation is that there is an
inductive structure to the definition whereby the ingredients for the
definition of algebraic $n$-stack involve only $n-1$-stacks and so are already
previously defined.

This definition came out of discussions with C. Walter in preparation for the
Trento school on algebraic stacks (September 1996). He made the remark that the
definition of algebraic stack made sense in any category where one has a
reasonable notion of smooth morphism, and suggested a general terminology of
``geometric stack'' for this notion.  One immediately realizes that the notion
of smooth morphism makes sense---notably---in the ``category'' of algebraic
stacks and therefore according to Walter's remark, one could define the notion
of geometric stack in the category of algebraic stacks.  This is the notion of
algebraic $2$-stack. It is an easy step to go from there to the general
inductive definition of algebraic $n$-stack.  Walter informs me that he had
also come upon the notion of algebraic $2$-stack at the same time (just
before the Trento school).

Now a note about terminology: I  have chosen to write the paper using Walter's
terminology ``geometric $n$-stack'' because this seems most closely to reflect
what is going on: the definition is made so that we can ``do geometry'' on the
$n$-stack, since in a rather strong sense it looks locally like a scheme. For
the purposes of the introduction, the terminology ``algebraic $n$-stack'' would
be better because this fits with Artin's terminology for $n=1$.  There is
another  place where the terminology ``algebraic'' would seem to be useful,
this is when we start to look at geometric stacks on the analytic site, which
we call ``analytic $n$-stacks''.  In fact one could interchange the
terminologies and in case of confusion one could even say ``algebraic-geometric
$n$-stack''.

In
\cite{RelativeLie} I proposed a notion of {\em presentable $n$-stack}
stable under homotopy fiber products and truncation.  One key part of the notion
of algebraic stack is the smoothness of the morphism $X\rightarrow T$ from a
scheme.  This is lost under truncation (e.g. the sheaf $\pi _0$ of an algebraic
stack may no longer be an algebraic stack); this indicates that the notion of
``geometric stack'' is something which combines together the various homotopy
groups in a fairly intricate way. In particular, the notion of presentable
$n$-stack is not the same as the notion of geometric $n$-stack (however a
geometric $n$-stack will be presentable).  This is a little bit analogous to the
difference between constructible and locally closed subsets in the theory of
schemes.

We will work over the site $\Xx$ of schemes of finite type over $Spec (k)$ with
the etale topology, and with the notion of smooth morphism. The definitions and
basic properties should also work for any site in which fiber products exist,
provided with a certain class of morphisms analogous to the smooth morphisms.
Rather than carrying this generalization through in the discussion, we leave it
to the reader.  Note that there are several examples which readily come to mind:
\newline
---the site of schemes of finite type with the etale topology and the class of
etale morphisms: this gives a notion of what might be called a ``Deligne-Mumford
$n$-stack'';
\newline
---the site of schemes of finite type with the fppf topology and
the class of flat morphisms: this gives a notion of  what might be called a
``flat-geometric $n$-stack'';
\newline
---the site of schemes of finite type with the qff topology and the class of
quasi-finite flat morphisms: this gives a notion of  what might be called
a ``qff-geometric $n$-stack''.

Whereas Artin proves \cite{ArtinInventiones} that
flat-geometric $1$-stacks are also smooth-geometric stacks (i.e. those defined
as we do here using smooth morphisms)---his proof is recounted in
\cite{LaumonMB}---it seems unlikely that the same would be true for $n$-stacks.
Artin's method also shows that qff-geometric $1$-stacks are Deligne-Mumford
stacks. However it looks like
Deligne-Mumford $n$-stacks are essentially just gerbs over Deligne-Mumford
$1$-stacks, while on the other hand in characteristic $p$ one could apply
Dold-Puppe (see below) to a complex of finite flat group schemes to get a fairly
non-trivial qff-algebraic $n$-stack. This seems to show that the
implication ``qff-geometric
$\Rightarrow $ Deligne-Mumford'' no longer holds for $n$-stacks.  This is why
it seems  unlikely that Artin's reasoning for the implication ``flat-geometric
$\Rightarrow$ smooth-geometric'' will work for $n$-stacks.

Here is the plan of the paper. In \S 1 we give the basic definitions of
geometric $n$-stack and smooth morphism of geometric $n$-stacks. In \S 2 we
give some basic properties which amount to having a good notion of geometric
morphism between $n$-stacks (which are themselves not necessarily geometric).
In \S 3 we briefly discuss some ways one could obtain geometric $n$-stacks by
glueing. In \S 4 we show that geometric $n$-stacks are presentable in the sense
of \cite{RelativeLie}. This is probably an important tool if one wants to do
any sort of Postnikov induction, since presentable $n$-stacks are closed under
the truncation processes which make up the Postnikov tower (whereas the
notion of geometric $n$-stack is not closed under truncation). In \S 5 we do a
preliminary version of what should be a more general Quillen theory. We treat
only the $1$-connected case, and then go on in the subsection ``Dold-Puppe'' to
treat the relative (i.e. over a base scheme or base $n$-stack) stable (in the
sense of homotopy theory) case  in a
different way.  It would be nice to have a unified version including a
reasonable notion of differential graded Lie algebra over an $n$-stack $R$
giving an algebraic approach to relatively $1$-connected $n$-stacks over $R$,
but this seems a bit far off in a technical sense.

In \S 6 we look at maps from a projective variety (or a smooth formal category)
into a geometric $n$-stack.  Here again it would be nice to have a fairly
general theory covering maps into any geometric $n$-stack but we can only say
something interesting in the easiest case, that of maps into {\em connected
very presentable $T$}, i.e. $n$-stacks with $\pi _0(T)=\ast$, $\pi _1(T)$ an
affine algebraic group scheme and $\pi _i(T)$ a vector space for $i\geq 2$.
(The terminology ``very presentable'' comes from \cite{RelativeLie}).  At the
end we speculate on how one might generalize to various other classes of $T$.

In \S 7 we briefly present an approach to defining the tangent stack to a
geometric $n$-stack. This is a generalization of certain results in the last
chapter of \cite{LaumonMB} although we don't refer to the cotangent complex.

In \S 8 we explain how to use geometric $n$-stacks as a framework for looking
at de Rham theory for higher nonabelian cohomology.  This is sort of a
synthesis of things that are in \cite{SantaCruz} and \cite{kobe}.

\bigskip

We assume known an adequate theory of $n$-stacks of groupoids over a site $\Xx$.
The main thing we will need is the notion of fiber product (which of course
means---as it always shall below---what one would often call the
``homotopy fiber
product'').

We work over an algebraically closed field $k$ of characteristic zero,
and sometimes directly over field $k=\cc$ of complex numbers.
Note however that the
definition makes sense over arbitrary base scheme and the ``Basic properties''
hold true there.

The term ``connected'' when applied to an $n$-stack means that the sheaf $\pi
_0$ is the final object $\ast$ (represented by $Spec (k)$). In the case of a
$0$-stack represented by  $Y$ this should not be confused with
connectedness of the scheme $Y$ which is a different question.

\numero{Definitions}

Let $\Xx$ be the site of schemes of finite type over $Spec (k)$ with the etale
topology. We will define the following notions: that an $n$-stack $T$ be {\em
geometric}; and that a morphism $T\rightarrow Z$ from a geometric $n$-stack to a
scheme be {\em smooth}.  We define these notions together by induction on $n$.
Start by saying that a $0$-stack (sheaf of sets) is {\em geometric} if it is
represented by an algebraic space. Say that a morphism $T\rightarrow Z$ from a
geometric $0$-stack to a scheme is {\em smooth} if it is smooth as a morphism of
algebraic spaces.

Now we give the inductive definitions:
say that an $n$-stack
$T$ is {\em geometric} if: \newline
GS1---for any schemes $X$ and $Y$ and morphisms $X\rightarrow T$,
$Y\rightarrow T$ the fiber product $X\times _TY$ (which is an $n-1$-stack)
is geometric using the inductive definition; and
\newline
GS2---there is a scheme $X$ and a morphism of $n$-stacks
$f:X\rightarrow T$ which is surjective on $\pi _0$ with the property that for
any scheme $Y$ and morphism $Y\rightarrow T$, the morphism
$$
X\times _TY\rightarrow Y
$$
from a geometric $n-1$ stack to the scheme $Y$ is smooth (using the inductive
definition).

If $T$ is a geometric $n$-stack we say that a morphism $T\rightarrow Y$ to a
scheme is {\em smooth} if for at least one
morphism $X\rightarrow T$ as in
GS2, the composed morphism $X\rightarrow Y$ is a smooth morphism of schemes.

This completes our inductive pair of definitions.

For $n=1$ we recover the notion of algebraic stack, and in fact our definition
is a straightforward generalization to $n$-stacks of Artin's definition of
algebraic stack.

The following lemma shows that the phrase ``for at least one'' in the definition
of smoothness can be replaced by the phrase ``for any''.

\begin{lemma}
\label{independence}
Suppose $T\rightarrow Y$ is a morphism from an $n$-stack to a scheme which
is smooth according to the previous definition,
and suppose that $U\rightarrow T$
is a morphism from a scheme such that for any scheme $Z\rightarrow T$, $U\times
_TZ\rightarrow Z$ is smooth (again according to the previous definition, as
a morphism from an $n-1$-stack to a scheme).
Then
$U\rightarrow Y$ is a smooth morphism of schemes.
\end{lemma}
{\em Proof:}
We prove this for
$n$-stacks by induction on $n$. Let $X\rightarrow T$ be the morphism as in GS2
such that $X\rightarrow Y$ is smooth.  Let $R= X\times _TU$. This is an
$n-1$-stack and the morphisms $R\rightarrow X$ and $R\rightarrow U$ are both
smooth as morphisms from $n-1$-stacks to schemes according to the above
definition. Let $W\rightarrow R$ be a surjection from a scheme as in property
GM2. By the present lemma applied inductively for $n-1$-stacks, the morphisms
$W\rightarrow X$ and $W\rightarrow U$ are smooth morphisms of schemes. But
the condition that $X\rightarrow Y$ is smooth implies that $W\rightarrow Y$ is
smooth, and then since $W\rightarrow U$ is smooth and surjective we get that
$U\rightarrow Y$ is smooth as desired. This argument doesn't work when $n=0$
but then $R$ is itself an algebraic space and the maps $R\rightarrow X$ (hence
$R\rightarrow Y$) and $R\rightarrow U$ are smooth maps of algebraic spaces;
this implies directly that $U\rightarrow Y$ is smooth.
\eop

The following lemma shows that these definitions don't change if we think of an
$n$-stack as an $n+1$-stack etc.

\begin{lemma}
\label{ntom}
Suppose $T$ is an $n$-stack which, when considered as an $m$-stack for some
$m\geq n$, is a geometric $m$-stack. Then $T$ is a geometric $n$-stack.
Similarly smoothness of a morphism $T\rightarrow Y$ to a scheme when $T$ is
considered as an $m$-stack implies smoothness when $T$ is considered as an
$n$-stack.
\end{lemma}
{\em Proof:} We prove this by induction on $n$ and then
$m$. The case $n=0$ and $m=0$ is clear. First treat the case $n=0$ and any $m$:
suppose $T$ is a sheaf of sets which is a geometric $m$-stack. There is a
morphism $X\rightarrow T$ with $X$ a scheme, such that if we set $R= X\times
_TX$ then $R$ is an $m-1$-stack smooth over $X$. However $R$ is again a sheaf of
sets so by the inductive statement for $n=0$ and $m-1$ we have that $R$ is an
algebraic space.  Furthermore the smoothness of the morphism $R\rightarrow X$
with $R$ considered as an $m-1$-stack implies smoothness with $R$ considered as
a $0$-stack.  In particular $R$ is an algebraic space with smooth maps to the
projections. Since the quotient of an algebraic space by a smooth equivalence
relation is again an algebraic space, we get that $T$ is an algebraic space
i.e. a geometric $0$-stack (and note by the way that $X\rightarrow T$ is a
smooth surjective map of algebraic spaces).  This proves the first statement for
$(0,m)$. For the second statement, suppose $T\rightarrow Y$ is a morphism to a
scheme $Y$ which is smooth as a morphism from an $m$-stack. Then choose the
smooth  surjective morphism
$X\rightarrow T$; as we have seen above this is a smooth morphism of algebraic
spaces. The definition of smoothness now is that $X\rightarrow Y$ is smooth.
But this implies that $T\rightarrow Y$ is smooth. This completes the
inductive step for $(0,m)$.

Now suppose we want to show the lemma for $(n,m)$ with $n\geq 1$ and suppose we
know it for all $(n', m')$ with $n'<n$ or $n'=n$ and $m'<m$.  Let $T$ be an
$n$-stack which is geometric considered as an $m$-stack. If $X,Y\rightarrow T$
are maps from schemes then $X\times _TY$ is an $n-1$-stack which is geometric
when considered as an $m-1$-stack; by the induction hypothesis it is geometric
when considered as an $n-1$-stack, which verifies GS1.  Choose a smooth
surjection $X\rightarrow T$ from a scheme as in property GS2 for $m$-stacks.
Suppose $Y\rightarrow T$ is any morphism from a scheme. Then $X\times _TY$ is an
$n-1$-stack with a map to $Y$ which is smooth considered as a map from
$m-1$-stacks. Again by the induction hypothesis it is smooth considered as a map
from an $n-1$-stack to a scheme, so we get GS2 for $n$-stacks. This completes
the proof that $T$ is geometric when considered as an $n$-stack.

Finally suppose $T\rightarrow Y$ is a morphism from an $n$-stack to a scheme
which is smooth considered as a morphism from an $m$-stack. Choose a surjection
$X\rightarrow T$ as in property GS2 for $m$-stacks; we have seen above that it
also satisfies the same property for $n$-stacks.  By definition of smoothness
of our original morphism from an $m$-stack, the morphism $X\rightarrow Y$ is
smooth as a morphism of schemes; this gives smoothness of $T\rightarrow Y$
considered as a morphism from an $n$-stack to a scheme. This finishes the
inductive proof of the lemma.
\eop

{\em Remarks:}
\newline
(1)\, We can equally well make a definition of {\em Deligne-Mumford $n$-stack}
by
replacing ``smooth'' in the previous definition with ''etale''. This gives an
$n$-stack whose homotopy group sheaves are finite...
\newline
(2)\, We could also make definitions of flat-geometric or qff-geometric
$n$-stack, by replacing the smoothness conditoin by flatness or quasifinite
flatness. If all of these notions are in question then we will denote the
standard one by ``smooth-geometric $n$-stack''. Not to be confused with ``smooth
geometric $n$-stack'' which means a smooth-geometric $n$-stack which is smooth!

We now complete our collection of basic definitions in some obvious ways.
We
say that a morphism of $n$-stacks $R\rightarrow T$ is {\em geometric}
if for any scheme $Y$ and map $Y\rightarrow T$ the fiber product $R\times _TY$
is a geometric $n$-stack.

We
say that a geometric morphism of $n$-stacks $R\rightarrow T$ is {\em smooth} if
for any scheme $Y$ and map $Y\rightarrow T$  the morphism $R\times
_TY\rightarrow Y$ is a smooth morphism in the sense of our inductive
definition.

\begin{lemma}
If $T\rightarrow Z$ is a morphism from an $n$-stack to a scheme then it is
smooth and geometric in the sense of the previous paragraph, if and only if $T$
is geometric and the morphism is smooth in the sense of our inductive
definition.
\end{lemma}
{\em Proof:}
Suppose that $T$ is geometric and the  morphism is smooth in the
sense of the previous paragraph.  Then applying that to the scheme $Z$ itself
we obtain that the morphism is smooth in the sense of the inductive definition.
On the other hand, suppose the morphism is smooth in the sense of the inductive
definition. Let $X\rightarrow T$ be a surjection as in GS2.  Thus $X\rightarrow
Z$ is smooth. For any scheme $Y\rightarrow Z$ we have that $X\times
_ZY\rightarrow T\times _ZY$ is surjective and smooth in the sense of the
previous
paragraph; but in this case (and using the direction we have proved above) this
is exactly the statement that it satisfies the conditions of GS2 for the stack
$T\times _ZY$.  On the other hand $X\times _ZY\rightarrow Y$ is smooth. This
implies (via the independence of the choice in the original definition of
smoothness which comes from \ref{independence})  that that $T\times
_ZY\rightarrow Y$ is smooth in the original sense. As this works for all $Y$, we
get that $T\rightarrow Z$ is smooth in the new sense.
\eop

\numero{Basic properties}

We assume that the propositions, lemmas and corollaries in this section are
known for $n-1$-stacks and we are proving them all in a gigantic induction for
$n$-stacks. On the other hand, in proving any statement we can use the
{\em previous} statements for the same $n$, too.

\begin{proposition}
\label{fiberprod}
If $R$, $S$ and $T$ are geometric $n$-stacks with morphisms
$R,T\rightarrow S$ then the fiber product
$R\times _ST$ is a geometric $n$-stack.
\end{proposition}
{\em Proof:}
Suppose $R$, $S$ and $T$ are geometric $n$-stacks with morphisms $R\rightarrow
S$ and $T\rightarrow S$.  Let $X\rightarrow R$, $Y\rightarrow S$ and
$Z\rightarrow T$ be smooth surjective morphisms from schemes.
Choose a smooth surjective morphism $W\rightarrow X\times _SZ$ from a scheme
(by axiom GS2 for $S$).

By base change of the morphism $Z\rightarrow T$, the morphism $X\times
_SZ\rightarrow X\times _ST$ is  a geometric smooth surjective morphism.
We first claim that the morphism $W\rightarrow X\times _ST$ is smooth.
To prove this, suppose $A\rightarrow X\times _ST$ is a morphism from a scheme.
Then $W\times _{X\times _ST}A\rightarrow A$ is the composition of
$$
W\times _{X\times _ST}A \rightarrow (X\times _SZ)\times _{X\times _ST}A
\rightarrow A.
$$
Both morphisms are geometric and smooth, and all three terms are $n-1$-stacks
(note that in the middle $(X\times _SZ)\times _{X\times _ST}A= A\times _TZ$).
By the composition result for $n-1$-stacks (Corollary \ref{geocomposition} below
with our global induction hypothesis) the composed morphism
$W\times _{X\times _ST}A\rightarrow A$ is smooth, and this for any $A$. Thus
$W\rightarrow X\times _ST$ is smooth.

Next we claim that the morphism $W\rightarrow R\times _ST$ is smooth. Again
suppose that $A\rightarrow R\times _ST$ is a morphism from a scheme. The two
morphisms
$$
W\times _{R\times _ST} A \rightarrow (X\times _ST)\times _{R\times _ST}A =
X\times _RA \rightarrow A
$$
are smooth and geometric by base change. Again this is a composition of
morphisms of $n-1$-stacks so by Corollary \ref{smoothcomposition2} and our
global
induction hypothesis the composition is smooth and geometric.  Finally the
morphism $W\rightarrow R\times _ST$ is the composition of three surjective
morphisms so it is surjective.  We obtain a morphism as in GS2 for $R\times
_ST$.

We turn to GS1.
Suppose $X\rightarrow R\times _ST$ and $Y\rightarrow R\times _ST$
are  morphisms from schemes. We would like to check that
$X\times _{R\times _ST}Y$ is a geometric $n-1$-stack. Note that calculating
$X\times _{R\times _ST}Y$ is basically the same thing as calculating in usual
homotopy theory the path space between two points $x$ and $y$ in a product of
fibrations $r\times _st$. From this point of view we see that
$$
X\times _{R\times _ST} Y =
(X\times _RY)\times _{X\times _SY}(X\times _TY).
$$
Note that the three components in the big fiber product on the right are
geometric $n-1$-stacks, so by our inductive hypothesis (i.e. assuming the
present proposition for $n-1$-stacks) we get that the right hand side is a
geometric $n-1$-stack, this gives the desired statement for GS1.
\eop

\begin{corollary}
If $R$ and $T$ are geometric $n$-stacks then any morphism between them is
geometric. In particular
an $n$-stack $T$ is geometric
if and only if the structural morphism $T\rightarrow \ast$ is geometric.
\end{corollary}
\eop

\begin{lemma}
\label{smoothcomposition1}
If $R\rightarrow S\rightarrow T$ are morphisms of geometric $n$-stacks and if
each morphism is smooth then the composition is smooth.
\end{lemma}
{\em Proof:}
We have already proved this for morphisms $X\rightarrow T \rightarrow Y$
where $X$ and $Y$ are schemes (see Lemma \ref{independence}). Suppose
$U\rightarrow T\rightarrow Y$ are smooth morphisms of geometric $n$-stacks
with
$Y$ a scheme. We prove that the composition is smooth, by induction on $n$ (we
already know it for $n=0$).  If
$Z\rightarrow T$ is a smooth surjective morphism from a scheme then the
morphism
$$
U\times _TZ \rightarrow Z
$$
is smooth by the definition of smoothness of $U\rightarrow T$. Also the map
$Z\rightarrow Y$ is smooth by definition of smoothness of $T\rightarrow Y$.
Choose a smooth surjection $V\rightarrow U\times _TZ$ from a scheme $V$ and note
that the map $V\rightarrow Z$ is smooth by definition, so (since these are
morphisms of schemes) the composition $V\rightarrow Y$ is smooth.
On the other
hand  $$ U\times _TZ \rightarrow U
$$
is smooth and surjective, by base change from
$Z\rightarrow T$.
We claim that the morphism $V\rightarrow U$ is smooth and surjective---actually
surjectivity is obvious.  To prove that it is smooth, let $W\rightarrow U$ be a
morphism from a scheme; then
$$
W\times _UV \rightarrow W\times _U(U\times _TZ) = W\times _TZ \rightarrow W
$$
is a composable pair of morphisms of $n-1$-stacks each of which is smooth
by base
change. By our induction hypothesis the composition is smooth.  This shows by
definition that $V\rightarrow U$ is smooth.
In particular the map $V\rightarrow U$ is admissible as  in GS2, and then we
can conclude that the map $U\rightarrow Y$ is smooth by the original definition
using $V$. This completes the proof in the current case.

Suppose finally that $U\rightarrow T \rightarrow R$ are smooth
morphisms of geometric $n$-stacks. Then for any scheme $Y$ the morphisms
$U\times _RY\rightarrow T\times _RY \rightarrow Y$ are smooth by base change;
thus from the case treated above their composition is smooth, and this
is the definition of smoothness of $U\rightarrow R$.
\eop

\begin{lemma}
\label{descendgeometric}
Suppose $S\rightarrow T$ is a geometric smooth surjective morphism of
$n$-stacks,
and suppose that $S$ is geometric. Then $T$ is geometric.
\end{lemma}
{\em Proof:}
We first show GS2. Let $W\rightarrow S$ be a smooth geometric surjection from a
scheme. We claim that the morphism $W\rightarrow T$ is surjective (easy),
geometric and smooth.  To show that it is geometric, suppose $Y\rightarrow T$
is a morphism from a scheme. Then since $S\rightarrow T$ is geometric we have
that $Y\times _TS$ is a geometric $n$-stack. On the other hand,
$$
Y\times _TW = (Y\times _TS)\times _SW,
$$
so by Proposition \ref{fiberprod} $Y\times _TW$ is geometric.  Finally to show
that $W\rightarrow T$ is smooth, note that
$$
Y\times _TW\rightarrow Y\times _TS \rightarrow Y
$$
is a  composable pair of smooth (by base change) morphisms of geometric
$n$-stacks, so by the previous lemma the composition is smooth. The morphism
$W\rightarrow T$ thus works for condition GS2.

To show GS1, suppose $X,Y\rightarrow T$ are morphisms from schemes. Then
$$
(X\times _TY)\times _TW = (X\times _TW)\times _W (Y\times _TW).
$$
The geometricity of the morphism $W\rightarrow T$ implies that $X\times _TW$
and $Y\times _TW$ are geometric, whereas of course $W$ is geometric. Thus by
Proposition \ref{fiberprod} we get that
$(X\times _TY)\times _TW$ is geometric. Now note that the morphism
$$
(X\times _TY)\times _TW \rightarrow X\times _TY
$$
of $n-1$-stacks is geometric, smooth and surjective (by base change of the same
properties for $W\rightarrow T$).  By the inductive version of the present
lemma for $n-1$ (noting that the lemma is automatically true for $n=0$) we
obtain that $X\times _TY$ is geometric. This is GS1.
\eop

\begin{corollary}
\label{localization}
Suppose $Y$ is a scheme and $T\rightarrow Y$ is a morphism from an $n$-stack. If
there is a smooth surjection $Y' \rightarrow Y$ such that $T':=Y'\times
_YT\rightarrow Y'$ is  geometric then the original morphism is geometric.
\end{corollary}
{\em Proof:}
The morphism $T'\rightarrow T$ is geometric, smooth and surjective (all by
base-change from the morphism $Y'\rightarrow Y$).  By \ref{descendgeometric},
the
fact that $T'$ is geometric implies that $T$ is geometric.
\eop

This corollary is particularly useful to do etale localization.  It implies
that the property of a morphism of $n$-stacks being geometric, is etale-local
over the base.

\begin{corollary}
\label{fibration}
Given a geometric morphism $R\rightarrow T$ of $n$-stacks such that $T$ is
geometric, then $R$ is geometric.
\end{corollary}
{\em Proof:}
Let $X\rightarrow T$ be the geometric smooth surjective morphism from a scheme
given by GS2 for $T$. By base change, $X\times _TR \rightarrow R$ is a geometric
smooth surjective morphism. However, by the geometricity of the morphism
$R\rightarrow T$ the fiber product $X\times _TR$ is geometric; thus by the
previous lemma, $R$ is geometric.
\eop

\begin{corollary}
\label{geocomposition}
The composition of two geometric  morphisms is
again geometric.
\end{corollary}
{\em Proof:}
Suppose $U\rightarrow T\rightarrow R$ are geometric morphisms, and suppose
$Y\rightarrow R$ is a morphism from a scheme. Then
$$
U\times _RY= U\times _T(T\times _RY).
$$
By hypothesis $T\times _RY$ is geometric.  On the other hand $U\times
_RY\rightarrow T\times _RY$ is geometric (since the property of being geometric
is obviously stable under base change).  By the previous Proposition
\ref{fiberprod} we get that $U\times _RY$ is geometric. Thus the morphism
$U\rightarrow R$ is geometric.
\eop

\begin{corollary}
\label{smoothcomposition2}
The composition of two  geometric smooth morphisms is
geometric and smooth.
\end{corollary}
{\em Proof:}
Suppose $R\rightarrow S \rightarrow T$ is a pair of geometric smooth morphisms.
Suppose  $Y\rightarrow T$ is a morphism from a scheme. Then
(noting by the previous corollary that $R\rightarrow T$ is geometric)
$R\times _TY$ and $S\times _TY$ are geometric. The composable pair
$$
R\times _TY \rightarrow S\times _TY \rightarrow Y
$$
of smooth morphisms now falls into the hypotheses of Lemma
\ref{smoothcomposition1} so the composition is smooth.  This implies that our
original composition was smooth.
\eop

In a relative setting we get:
\begin{corollary}
Suppose $U\stackrel{a}{\rightarrow}T\stackrel{b}{\rightarrow}R$ is a composable
pair of morphisms of $n$-stacks. If $a$ is geometric, smooth and surjective
and $ba$ is geometric (resp. geometric and smooth) then $b$ is geometric (resp.
geometric and smooth).
\end{corollary}
{\em Proof:}
Suppose $Y\rightarrow R$ is a morphism from a scheme. Then
$$
Y\times _RU = (Y\times _RT)\times _TU.
$$
The map $Y\times _RU\rightarrow Y\times _RT$ is geometric, smooth and
surjective (since those properties are obviously---from the form of their
definitions---invariant under base change). The fact that $ba$ is geometric
implies that $Y\times _RU$ is geometric, which by the previous lemma implies
that $Y\times _RT$ is geometric. Suppose furthermore that $ba$ is smooth.
Choose a smooth surjection $W\rightarrow Y\times _RT$ from a scheme. Then
the morphism
$$
W\times _{Y\times _RT} (Y\times _RU)\rightarrow Y\times _RU
$$
is smooth by basechange and the morphism $Y\times _RU\rightarrow Y$ is smooth
by hypothesis. Thus $W\times _{Y\times _RT} (Y\times _RU)\rightarrow Y$
is smooth.  Choosing a smooth surjection from a scheme
$$
V\rightarrow W\times _{Y\times _RT} (Y\times _RU)
$$
we get that $V\rightarrow Y$ is a smooth morphism of schemes.
On the other hand, the morphism
$$
W\times _{Y\times _RT} (Y\times _RU)\rightarrow W
$$
is smooth and surjective, so $V\rightarrow W$ is smooth and surjective.
Therefore $W\rightarrow Y$ is smooth. This proves that if $ba$ is smooth then
$b$ is smooth.
\eop

{\em Examples:}
Proposition \ref{fibration} allows us to construct many
examples. The main examples we shall look at below are the {\em connected
presentable $n$-stacks}.  These are connected $n$-stacks $T$ with (choosing a
basepoint $t\in T(Spec (\cc ))$ $\pi _i (T, t)$ represented by
group schemes of finite type.  We apply \ref{fibration} inductively to
show that such a $T$ is geometric. Let $T\rightarrow \tau _{\leq n-1}T$ be the
truncation morphism.  The fiber over a morphism $Y\rightarrow \tau _{\leq
n-1}T$ is (locally in the etale topology of $Y$ where there exists a
section---this is good enough by \ref{localization}) isomorphic to $K(G/Y, n)$
for a smooth group scheme of finite type $G$ over $Y$.  Using the following
lemma, by induction $T$ is geometric.

\begin{lemma}
\label{eilenbergExample}
Fix $n$, suppose $Y$ is a scheme and suppose $G$ is a smooth group scheme over
$Y$. If $n\geq 2$ require $G$ to be abelian.  Then $K(G/Y, n)$ is a geometric
$n$-stack and the morphism $K(G/Y,n)\rightarrow Y$ is smooth..
\end{lemma}
{\em Proof:}
We prove this by induction on $n$. For $n=0$ we simply have $K(G/Y,0)=G$ which
is a scheme and hence geometric---also note that by hypothesis it is smooth over
$Y$.  Now for any $n$, consider the basepoint section $Y\rightarrow K(G/Y,n)$.
We claim that this is a smooth geometric map.  If $Z\rightarrow K(G/Y,n)$ is
any morphism then it corresponds to a map $Z\rightarrow Y$ and a class in
$H^n(Z,G|_Z)$. Since we are working with the etale topology, by definition this
class vanishes on an etale surjection $Z'\rightarrow Z$ and for our claim it
suffices to show that $Y\times _{K(G/Y,n)}Z'$ is smooth and geometric over
$Z'$.  Thus we may assume that our map $Z'\rightarrow K(G/Y,n)$ factors through
the basepoint section $Y\rightarrow K(G/Y,n)$. In particular it suffices to
prove that $Y\times _{K(G/Y,n)}Y\rightarrow Y$ is smooth and geometric. But $$
Y\times _{K(G/Y,n)}Y= K(G/Y,n-1)
$$
so by our induction hypothesis this is geometric and smooth over $Y$. This
shows that $K(G/Y,n)$ is geometric and furthermore the basepoint section is a
choice of map as in GS2.  Now the composed map $Y\rightarrow K(G/Y,n)\rightarrow
Y$ is the identity, in particular smooth, so by definition $K(G/Y,n)\rightarrow
Y$ is smooth.
\eop

Note that stability under fiber products (Proposition \ref{fiberprod}) implies
that if $T$ is a geometric $n$-stack then $Hom (K, T)$ is geometric for any
finite CW complex $K$. See (\cite{kobe} Corollary 5.6) for the details of the
argument---which was in the context of presentable $n$-stacks but the argument
given there only depends on stability of our class of $n$-stacks under fiber
product.  We can apply this in particular to the geometric $n$-stacks
constructed just above, to obtain some non-connected examples.

If $T= BG$ for an algebraic group $G$ and $K$ is connected with basepoint
$k$ then $Hom (K, T)$ is the moduli stack of representations $\pi
_1(K,k)\rightarrow G$ up to conjugacy.

\numero{Locally geometric $n$-stacks}

The theory we have described up till now concerns objects {\em of finite type}
since we have assumed that the scheme $X$ surjecting to our $n$-stack $T$ is of
finite type.  We can obtain a definition of ``locally geometric'' by relaxing
this to the condition that $X$ be locally of finite type (or equivalently that
$X$ be a disjoint union of schemes of finite type). To be precise we say that an
$n$-stack $T$ is {\em locally geometric} if there exists a sheaf
which is a disjoint union of schemes of finite type, with a morphism
$$
 \varphi : X=\coprod
X_i \rightarrow T
$$
such that $\varphi$ is smooth and geometric.

Note that if $X$ and $Y$ are schemes of finite type mapping to $T$ we still
have that $X\times _TY$ is geometric (GS1).

All of the previous results about fibrations, fiber products, and so on still
hold for locally geometric $n$-stacks.

One might want also to relax the definition even further by only requiring
that $X\times _TY$ be itself locally geometric (and so on) even for schemes of
finite type.  We can obtain a notion that we call {\em slightly
geometric} by replacing ``scheme of finite type'' by ``scheme locally
of finite type'' everywhere in the preceeding definitions.  This notion may be
useful in the sense that a lot more $n$-stacks will be ``slightly geometric''.
However it seems to remove us somewhat from the realm where geometric reasoning
will work very well.

\numero{Glueing}

We say that a morphism $U\rightarrow T$ of geometric stacks is a {\em
Zariski open subset} (resp. {\em etale open subset}) if for every scheme $Z$
and $Z\rightarrow T$ the fiber product $Z\times _TU$ is a Zariski open subset
of $Z$ (resp. an algebraic space with etale map to $Z$).

If we have two
geometric $n$-stacks $U$ and $V$ and a geometric  $n$-stack $W$ with morphisms
$W\rightarrow U$ and $W\rightarrow V$ each of which is a Zariski open subset,
then we can glue $U$ and $V$ together along $W$ to get a geometric $n$-stack
$T$ with Zariski open subsets $U\rightarrow T$ and $V\rightarrow T$ whose
intersection is $W$. If one wants to glue several open sets it has to be done
one at a time (this way we avoid having to talk about higher cocycles).

As a more general result we have the following.  Suppose $\Phi$ is a functor
from the simplicial category $\Delta$ to the category of $n$-stacks (say a
strict functor to the category of simplicial presheaves, for example).
Suppose that each $\Phi _k$ is a geometric $n$-stack, and suppose that
the two morphisms $\Phi _1 \rightarrow \Phi _0$ are smooth. Suppose furthermore
that $\Phi $ satisfies the Segal condition that
$$
\Phi _k \rightarrow \Phi _1\times _{\Phi _0} \ldots \times _{\Phi _0}\Phi _1
$$
is an equivalence (i.e. Illusie weak equivalence of simplicial presheaves).
Finally suppose that for any element of $\Phi _1(X)$ there is, up to
localization over $X$, an ``inverse'' (for the multiplication on $\Phi _1$ that
comes from Segal's condition as in \cite{Segal}) up to homotopy.
Let $T$
be the realization over the simplicial variable, into a presheaf of spaces (i.e.
we obtain a bisimplicial presheaf, take the diagonal).

\begin{proposition}
In the above situation, $T$ is a geometric $n+1$-stack.
\end{proposition}
{\em  Proof:}
There is a surjective map $\Phi _0 \rightarrow T$ and we have by definition that
$\Phi _0 \times _T\Phi _0 = \Phi _1$.  From this one can see that $T$ is
geometric.
\eop

As an example of how to apply the above result, suppose $U$ is a geometric
$n$-stack and suppose we have a geometric $n$-stack $R$ with
$R\rightarrow U\times U$. Suppose furthermore that we have a multiplication
$R\times _{p_2, U, p_1}R\rightarrow R$ which is associative and such that
inverses exist up to homotopy.  Then we can set $\Phi _k = R\times _U\ldots
\times _UR$ with $\Phi _0 = U$. We are in the above situation, so we obtain the
geometric $n$-stack $T$. We call this the {\em $n$-stack associated to the
descent data $(U, R)$}.

The original result about glueing over Zariski open
subsets can be interpreted in this way.

The simplicial version of this descent
with any $\Phi$ satisfying Segal's condition is a way to avoid having to talk
about strict associativity of the composition on $R$.

\numero{Presentability}
Recall first of all that the category of {\em vector sheaves} over a scheme $Y$
is the smallest abelian category of abelian sheaves on the big site of
schemes over $Y$ containing the structure sheaf (these were called
``$U$-coherent
sheaves'' by Hirschowitz in \cite{Hirschowitz}, who was the first to define
them). A vector sheaf may be presented as the kernel of a sequence of $3$
coherent sheaves which is otherwise exact on the big site; or dually as the
cokernel of an otherwise-exact sequence of $3$ {\em vector schemes} (i.e. duals
of coherent sheaves).  The nicest thing about the category of vector sheaves is
that duality is involutive.

Recall that we have defined in \cite{RelativeLie} a notion of {\em presentable
group sheaf} over any base scheme $Y$.  We will  not repeat the definition
here, but just remark (so as to give a rough idea of what is going on) that if
$G$ is a presentable group sheaf over $Y$ then it admits a Lie algebra object
$Lie (G)$ which is a vector sheaf with bilinear Lie bracket operation
(satisfying Jacobi).

In \cite{RelativeLie} a definition was then made of {\em presentable $n$-stack};
this involves a certain condition on $\pi _0$ (for which  we refer to
\cite{RelativeLie})
and the condition that the higher homotopy group sheaves (over any base scheme)
be presentable group sheaves.

For our purposes we shall often be interested in the slightly more restrictive
notion of {\em very presentable $n$-stack}.  An $n$-stack $T$ is defined (in
\cite{RelativeLie}) to be very presentable if it is presentable, and if
furthermore:
\newline
(1)\, for $i\geq 2$ and for any scheme $Y$ and $t\in T(Y)$ we have
that $\pi _i (T|_{\Xx /Y}, t)$ is a vector sheaf over $Y$; and
\newline
(2)\, for any artinian scheme $Y$ and $t\in T(Y)$ the group of sections
$\pi _1(T|_{\Xx /Y},t)(Y)$ (which is naturally an algebraic group scheme over
$Spec (k)$) is affine.

For our purposes here we will mostly stick to the case of connected $n$-stacks
in the coefficients. Thus we review what the above definitions mean for $T$
connected (i.e. $\pi _0(T)=\ast $). Assume that $k$ is algebraically closed
(otherwise one has to take what is said below possibly with some Galois
twisting). In the connected case there is essentially a unique basepoint $t\in
T(Spec (k))$.  A group sheaf over $Spec (k)$ is presentable if and only if it is
an algebraic group scheme (\cite{RelativeLie}), so $T$ is presentable if
and only
if $\pi _i(T,t)$ are represented by algebraic group schemes.  Note that a
vector sheaf over $Spec (k)$ is just a vector space, so $T$ is very presentable
if and only if the $\pi _i (T,t)$ are vector spaces for $i\geq 2$ and $\pi
_1(T,t)$ is an affine algebraic group scheme (which can of course act on the
$\pi _i$ by a representation which---because we work over the big
site---is automatically algebraic).

\begin{proposition}
\label{presentable}
If $T$ is a geometric $n$-stack on $\Xx$ then $T$ is presentable in the sense
of \cite{RelativeLie}.
\end{proposition}
{\em Proof:}
Suppose $X\rightarrow R$ is a smooth morphism from a scheme $X$ to a geometric
$n$-stack $R$.  Note that the morphism
$R\rightarrow \pi _0(R)$ satisfies the lifting properties $Lift _n(Y, Y_i)$,
since by localizing in the etale topology we get rid of any cohomological
obstructions to lifting coming from the higher homotopy groups. On the other
hand the morphism $X\rightarrow R$ being smooth, it satisfies the lifting
properties (for example one can say that the map $X\times _RY\rightarrow Y$
is smooth and admits a smooth surjection from a scheme smooth over $Y$;
with this one gets the lifting properties, recalling of course that a smooth
morphism between schemes is vertical. Thus we get that
$X\rightarrow \pi _0(R)$ is vertical.

Now suppose $T$ is geometric and choose a smooth surjection $u:X\rightarrow T$.
We get from above that $X\rightarrow \pi _0(T)$ is vertical.  Note that
$$
X\times _{\pi _0(T)}X = im (X\times _TX \rightarrow X\times X).
$$
Let $G$ denote the group sheaf $\pi _1(T|_{\Xx /X}, u)$ over $X$.
We have that $G$ acts freely on $\pi _0(X\times _TX)$ (relative to the first
projection $X\times _TX\rightarrow X$) and the quotient is the image
$X\times _{\pi _0(T)}X$. Thus, locally over schemes mapping into
the target, the morphism
$$
\pi _0(X\times _TX) \rightarrow X\times _{\pi _0(T)}X
$$
is the same as the morphism
$$
G\times _X(X\times _{\pi _0(T)}X)\rightarrow X\times _{\pi _0(T)}X
$$
obtained by base-change. Since $G\rightarrow X$ is a group sheaf it is an
$X$-vertical morphism (\cite{RelativeLie} Theorem 2.2 (7)), therefore its base
change is again an $X$-vertical morphism. Since verticality is local over
schemes mapping into the target, we get that
$$
\pi _0(X\times _TX) \rightarrow X\times _{\pi _0(T)}X
$$
is an $X$-vertical morphism. On the other hand by the definition that $T$ is
geometric we obtain a smooth
surjection $R\rightarrow X\times _TX$ from a scheme
$R$, and by the previous discussion this gives a $Spec(\cc )$-vertical
surjection
$$
R\rightarrow \pi _0(X\times _TX).
$$
Composing we get the $X$-vertical surjection $R\rightarrow X\times _{\pi
_0(T)}X$. We have now proven that $\pi _0(T)$ is $P3\frac{1}{2}$ in the
terminology of \cite{RelativeLie}.

Suppose now that $v:Y\rightarrow T$ is a point. Let $T':= Y\times _TY$.
Then $\pi _0(T')= \pi _1(T|_{\Xx /Y}, v)$ is the group sheaf we are interested
in looking at over $Y$. We will show that it is presentable.
Note that $T'$ is
geometric; we apply the same argument as above, choosing a smooth surjection
$X\rightarrow T'$. Recall that this gives a $Spec (\cc )$-vertical (and hence
$Y$-vertical) surjection $X\rightarrow \pi _0(T')$. Choose a smooth surjection
$R\rightarrow X\times _{T'}X$. In the previous proof the group sheaf denoted $G$
on $X$ is actually pulled back from a group sheaf $\pi _2(T|_{\Xx /Y}, v)$ on
$Y$.  Therefore the morphism
$$
\pi _0(X\times _{T'}X)\rightarrow X\times _{\pi _0(T')}X
$$
is a quotient by a group sheaf over $Y$, in particular it is $Y$-vertical.
As usual the morphism $R\rightarrow \pi _0(X\times _{T'}X)$ is $Spec (\cc
)$-vertical so in particular $Y$-vertical. We obtain a $Y$-vertical surjection
$$
R\rightarrow X\times _{\pi _0(T')}X.
$$
This finishes the proof that $\pi _1(T|_{\Xx /Y}, v)$
satisfies property $P4$ (and since it
is a group sheaf, $P5$ i.e. presentable) with respect to $Y$.

Now note that $\pi _i(T|_{\Xx /Y}, v) = \pi _{i-1}(T'|_{\Xx /Y}, d)$
where $d: Y \rightarrow T' := Y\times _TY$ is the diagonal morphism.
Hence (as $T'$ is itself geometric) we obtain by induction that all of the
$\pi _i(T|_{\Xx /Y}, v)$ are presentable group sheaves over $Y$.
This shows that $T$ is presentable in the terminology of \cite{RelativeLie}.
\eop

Note that presentability in \cite{RelativeLie} is  a slightly stronger condition
than the condition of presentability as it is referred to in \cite{kobe} so all
of the results stated in \cite{kobe} hold here; and of course all of the
results of \cite{RelativeLie} concerning presentable $n$-stacks hold for
geometric $n$-stacks. The example given below which shows that the class of
geometric $n$-stacks is not closed under truncation, implies that the class of
presentable $n$-stacks is strictly bigger than the class of geometric ones,
since the class of presentable $n$-stacks is closed under truncation
\cite{RelativeLie}.

The results announced (some with sketches of proofs) in
\cite{kobe} for presentable $n$-stacks hold for geometric $n$-stacks.
Similarly the basic results of \cite{RelativeLie} hold for geometric $n$-stacks.
For example, if $T$ is a geometric $n$-stack and $f:Y\rightarrow T$ is a
morphism
from a scheme then $\pi _i (T|_{\Xx /Y}, f)$ is a presentable group sheaf, so
it has a Lie algebra object $Lie\, \pi _i (T|_{\Xx /Y}, f)$ which is a {\em
vector sheaf} (or ``$U$-coherent sheaf'' in the terminology of
\cite{Hirschowitz}) with Lie bracket operation.

{\em Remark:}  By Proposition \ref{presentable}, the condition of being
geometric
is stronger than the condition of being presentable given in
\cite{RelativeLie}.
Note from the example given below showing that geometricity is not compatible
with truncation (whereas by definition presentability is compatible with
truncation),
the condition of being geometric is {\em strictly} stronger than the condition
of being presentable.

Of course in the connected case, presentability and geometricity are the same
thing.

\begin{corollary}
A connected $n$-stack $T$ is geometric if and only if the $\pi _i(T,t)$ are
group schemes of finite type for all $i$.
\end{corollary}
{\em Proof:}
We show in \cite{RelativeLie} that presentable groups over $Spec (k)$ are just
group schemes of finite type.  Together with the previous result this shows
that if $T$ is connected and geometric then the $\pi _i(T,t)$ are
group schemes of finite type for all $i$. On the other hand, if $\pi
_0(T)=\ast$ and the  $\pi _i(T,t)$ are
group schemes of finite type for all $i$ then by the Postnikov decomposition
of $T$
and using \ref{fibration}, we conclude that $T$ is geometric (note that
for a group scheme of finite type $G$,
$K(G,n)$ is geometric).
\eop

\numero{Quillen theory}

Quillen in \cite{Quillen} associates to every $1$-connected rational
space $U$ a {\em differential graded Lie algebra (DGL)} $L_{\cdot} = \lambda
(U)$: a DGL is a graded Lie algebra
(over $\qq$ for our purposes) $L_{\cdot} = \bigoplus _{p\geq 1}L_p$
(with all elements of strictly positive degree)
with differential $\partial : L_p \rightarrow L_{p-1}$ compatible in the usual
(graded) way with the Lie bracket.  Note our conventions that the indexing is
downstairs, by positive numbers and the differential has degree $-1$. The
homology groups of $\lambda (U)$ are the homotopy groups of $U$ (shifted by one
degree).

This construction gives an equivalence between the homotopy theory of DGL's and
that of rational spaces.  Let $L_{\cdot} \mapsto |L_{\cdot}| $ denote the
construction going in the other direction.  We shall assume for our purposes
that there exists such a realization functor from the category of DGL's to
the category of $1$-connected spaces, compatible with finite direct products.

Let $DGL_{\cc , n}$ denote the category of $n$-truncated $\cc$-DGL's of finite
type (i.e. with homology groups which are finite dimensional vector
spaces, vanishing in degree $\geq n$) and free as graded Lie algebras.

We define a realization functor $\rho ^{\rm pre}$ from $DGL_{\cc ,n}$ to the
category of presheaves of spaces over $\Xx$.  If $L_{\cdot}\in DGL_{\cc
,n}$ then
for any $Y\in \Xx$ let
$$
\rho ^{\rm pre}(L_{\cdot})(Y):= | L_{\cdot} \otimes _{\cc}\Oo (Y) |.
$$
Then let $\rho (L_{\cdot})$ be the $n$-stack associated to the presheaf of
spaces
$\rho ^{\rm pre}(L_{\cdot})$.  This construction is functorial and compatible
with direct products (because we have assumed the same thing about the
realization functor $|L_{\cdot}|$).

Note that $\pi _0^{\rm pre}(\rho ^{\rm pre}(L_{\cdot}))=\ast$
and in fact we can choose a basepoint $x$ in $\rho ^{\rm pre}(L_{\cdot})(Spec
(\cc ))$.
We have
$$
\pi _i ^{\rm pre}(\rho ^{\rm pre}(L_{\cdot}), x)= H_{i-1}(L_{\cdot})
$$
(in other words the presheaf on the left is represented by the vector space on
the right).  This gives the same result on the level of associated stacks and
sheaves:
$$
\pi _i (\rho (L_{\cdot}), x)= H_{i-1}(L_{\cdot}).
$$
In particular note that a morphism of DGL's induces an equivalence of
$n$-stacks if and only if it is a quasiisomorphism.  Note also
that $\rho (L_{\cdot})$ is a $1$-connected $n$-stack whose higher homotopy
groups are complex vector spaces, thus it is a very presentable
geometric $n$-stack.

\begin{theorem}
The above construction gives an equivalence between the homotopy
category $ho\, DGL_{\cc , n}$ and the homotopy category of $1$-connected very
presentable $n$-stacks.
\end{theorem}
{\em Proof:}
Let $(L,M)$ denote the set of homotopy classes of maps from $L$ to $M$ (either
in the world of DGL's or in the world of $n$-stacks on $\Xx$). Note that if $L$
and $M$ are DGL's then $L$ should be free as a Lie algebra (otherwise we have
to replace it by a quasiisomorphic free one). We prove that the map
$$
(L,M)\rightarrow (\rho (L), \rho (M))
$$
is an isomorphism.  First we show this for the case where $L=V[n-1]$ and
$M=U[n-1]$ are finite dimensional vector spaces in degrees $n-1$ and $m-1$. In
this case (where unfortunately $L$ isn't free so has to be replaced by a free
DGL) we have
$$
(V[n-1], U[m-1])= Hom (Sym ^{m/n}(V), U)
$$
where the symmetric product is in the graded sense (i.e. alternating or
symmetric according to parity) and defined as zero when the exponent is not
integral.  Note that $\rho (V[n-1])= K(V, n)$ and $\rho (U[m-1])= K(U,m)$.
The Breen calculations in characteristic zero (easier than the case treated in
\cite{BreenIHES}) show that
$$
(K(V,n), K(U,m))= Hom (Sym ^{m/n}(V), U)
$$
so our claim holds in this case.

Next we treat the case of arbitrary $L$ but $M= U[m-1]$ is again a vector space
in degree $m-1$.  In this case we are calculating the cohomology of $L$ or
$\rho (L)$ in degree $m$ with coefficients in $U$.  Using a Postnikov
decomposition of $L$ and the appropriate analogues of the Leray spectral
sequence on both sides we see that our functor induces an isomorphism on these
cohomology groups.

Finally we get to the case of arbitrary $L$ and arbitrary $M$. We proceed by
induction on the truncation level $m$ of $M$. Let $M'=\tau _{\leq m-1}M$ be the
truncation  (coskeleton) with the natural morphism $M\rightarrow M'$. The fiber
is of the form $U[m-1]$ (we index our truncation by the usual homotopy groups
rather than the homology groups of the DGL's which are shifted by $1$). Note
that $\rho (M')= \tau _{\leq m-1}\rho (M)$ (since the construction $\rho$ is
compatible with homotopy groups so it is compatible with the truncation
operations). The fibration $M\rightarrow M'$ is classified by a map
$f:M'\rightarrow U[m]$ and the fibration $\rho (M)\rightarrow \rho (M')$
by the corresponding map $\rho (f): \rho (M')\rightarrow K(U, m+1)$.
The image of
$$
(L,M)\rightarrow (L, M')
$$
consists of the morphisms whose composition into $U[m]$ is homotopic to the
trivial morphism $L\rightarrow U[m]$. Similarly the image of
$$
(\rho (L),\rho (M))\rightarrow (\rho (L),\rho (M'))
$$
is the morphisms whose composition into $K(U,m+1)$ is homotopic to the
trivial morphism.  By our inductive hypothesis
$$
\rho : (L,M')\rightarrow (\rho (L), \rho (M'))
$$
is an isomorphism. The functor $\rho$ is an
isomorphism on the images, because we know the statement for targets $U[m]$.
Suppose we are given a map $a:L\rightarrow M'$ which is in the image.
The inverse
image of this homotopy class in $(L,M)$ is the quotient of the set of liftings
of $a$ by the action of the group of self-homotopies of the map $a$. The set of
liftings is a principal homogeneous space under $(L, U[m-1])$.

Similarly the inverse image of the homotopy class of $\rho (a)$ in $(\rho (L),
\rho (M))$ is the quotient of the set of liftings of $\rho (a)$ by the group of
self-homotopies of $\rho (a)$. Again the set of liftings is a principal
homogeneous space under $(\rho (L), K(U,m))$.

The actions in the principal homogeneous spaces come from maps
$$
U[m-1]\times M\rightarrow M
$$
over $M'$ and
$$
K(U, m)\times \rho (M)\rightarrow \rho (M)
$$
over $\rho (M')$, the second of which is the image under $\rho$ of the first.

Since $\rho : (L, U[m-1])\cong (\rho (L), K(U,m))$, we will get that $\rho$
gives an isomorphism of the fibers if we can show that the images of the actions
of the groups of self-homotopy equivalences are the same. Notice that since
these actions are on principal homogeneous spaces they factor through the
abelianizations of the groups of self-homotopy equivalences.

In general if $A$ and $B$ are spaces then $(A\times S^1, B)$ is the disjoint
union over $(A,B)$ of the sets of conjugacy classes of the groups of
self-homotopies  of the maps from $A$ to $B$. On the other hand a map of groups
$G\rightarrow G'$ which induces an isomorphism on sets of conjugacy classes is
surjective on the level of abelianizations. Thus if we know that a certain
functor gives an isomorphism on $(A,B)$ and on $(A\times S^1, B)$ then it is a
surjection on the abelianizations of the groups of self-homotopies of the maps.

Applying this principle in the above situation, and noting that we know
by our induction hypothesis that $\rho$ induces isomorphisms on $(L, M')$ and
$(L\times \lambda (S^1)\otimes _{\zz}k, M')$, we find that $\rho$ induces a
surjection from the abelianization of the group of self-homotopies of the map
$a:L\rightarrow M'$ to the abelianization of the group of self-homotopies of
$\rho (a)$.  This finally allows us to conclude that $\rho$ induces an
isomorphism from the inverse image of the class of $a$ in $(L,M)$ to the
inverse image of the class of $\rho (a)$ in $(\rho (L), \rho (M))$. We have
completed our proof that
$$
\rho : (L,M)\cong (\rho (L), \rho (M)).
$$

In order to obtain that $\rho$ induces an isomorphism on homotopy categories we
just have to see that any $1$-connected very presentable $n$-stack $T$ is of
the form $\rho (L)$. We show this by induction on the truncation level.
Put  $T'=\tau _{\leq n-1}T$. By the induction hypothesis there is a DGL $L'$
with $\rho (L')\cong T'$ (and we may actually write $\rho (L')=T'$).  Now the
fibration $T\rightarrow T'$ is classified by a map $f:T'\rightarrow K(V,n+1)$.
From the above proof this map comes from a map $b:L'\rightarrow V[n]$, that is
$f=\rho (b)$.   In turn this map classifies a new
DGL $L$ over $L'$. The fibration $\rho (L)\rightarrow \rho (L')=T'$ is
classified by the map $\rho (b)=f$ so $\rho (L)\cong T$.
\eop

\subnumero{Dold-Puppe}

Eventually it would be nice to have a relative version of the previous theory,
over any $n$-stack $R$. The main problem in trying to do this is to
have the right notion of complex of sheaves over an $n$-stack $R$.  Instead of
trying to do this, we will simply use the notion of {\em spectrum over $R$} (to
be precise I will use the word ``spectrum'' for what is usually called an
``$\Omega$-spectrum''. For our purposes we are only interested in spectra with
homotopy groups which are rational and vanish outside of a bounded interval.  In
absolute terms such a spectrum is equivalent to a complex of rational vector
spaces, so in the relative case over a presheaf of spaces $R$ this gives a
generalized notion of complex over $R$.

For our spectra with homotopy groups nonvanishing only in a bounded region, we
can replace the complicated general theory by the simple consideration of
supposing that we are in the stable range.  Thus we fix numbers
$N, M$ with $M$ bigger than the length of any complex we want to consider and
$N$ bigger than $2M+2$.  For example if we are only interested in dealing with
$n$-stacks then we cannot be interested in complexes of length  bigger than $n$
so we could take $M>n$.

An {\em spectrum} (in our setting) is then simply  an $N$-truncated
rational space with a basepoint, and which is $N-M-1$-connected. More
generally if $R$ is an $n$-stack with $n \leq N$ then a {\em spectrum over
$R$} is just an $N$-stack $S$ with morphism $p:S\rightarrow R$ and section
denoted $\xi : R\rightarrow S$ such that  $S$ is rational and $N-M-1$-connected
relative to $R$.  A morphism of spectra is a morphism of spaces (preserving the
basepoint).

Suppose $S$ is a spectrum; we define the {\em complex associated to $S$}
by setting  $\gamma (S)^i$ to be the singular $N-i$-chains on $S$. The
differential $d: \gamma (S)^i\rightarrow \gamma (S)^{i+1}$ is the same as the
boundary map on chains (which switches direction because of the change in
indexing). Note that we have normalize things so that the complex starts in
degree $0$.  The homotopy theory of spectra is the same as that of complexes of
rational vector spaces indexed in degrees $\geq 0$,
with cohomology nonvanishing only in degrees $\leq M$.

If $C^{\cdot}$ is a complex as above then let $\sigma (C^{\cdot})$ denote the
corresponding spectrum.

This can be generalized to the case where the base is a $0$-stack. If $Y$ is
a $0$-stack (notably for example a scheme) and if $S$ is a spectrum over $Y$
then we obtain a complex of presheaves of rational vector spaces
$\gamma (S/Y)$ over $Y$. Conversely if $C^{\cdot}$ is a complex of presheaves
of rational vector spaces over $Y$ then we obtain a spectrum denoted $\sigma
(C^{\cdot}/Y)$. These constructions are an equivalence in homotopy theories,
where the weak equivalence between complexes means quasiisomorphism (i.e.
morphisms inducing isomorphisms on associated cohomology sheaves).

If $S$ is a spectrum and $n \leq N$ then we can define the {\em realization}
$\kappa (S, n)$ to be the $N-n$-th loop space $\Omega ^{N-n}S$ (the loops are
taken based at the given basepoint).  Similarly if $S$ is a spectrum over an
$n'$-stack $R$ then we obtain the {\em realization} $\kappa
(S/R,n)\rightarrow R$
as the $N-n$-th relative loop space based at the given section $\xi$.

Taken together we obtain the following construction: if $C^{\cdot}$ is a
complex of vector spaces then $\kappa (\sigma (C^{\cdot}), n)$ is an $n$-stack.
If $C^{\cdot}$ is a complex of presheaves of rational vector spaces over a
$0$-stack (presheaf of sets) $Y$ then $\kappa (\sigma (C^{\cdot}/Y)/Y, n)$
is an $n$-stack over $Y$.  These constructions are what is known as {\em
Dold-Puppe}. They are compatible with the usual Eilenberg-MacLane constructions:
if $V$ is a presheaf of rational vector spaces over $Y$ considered as a complex
in degree $0$ then $$
\kappa (\sigma (V/Y)/Y, n)= K(V/Y, n).
$$
The basic idea behind our notational system is that we think of spectra over
$R$ as being complexes of rational presheaves over $R$ starting in degree $0$.
The operation $\kappa (S/R, n)$ is the {\em Dold-Puppe} realization from a
``complex'' to a space relative to $R$.

We can do higher direct images in this context. If $f:R\rightarrow T$
is a morphism of $n$-stacks and if $S$ is a spectrum over $R$ then define
$f_{\ast}(S)$ to be the $N$-stack $\Gamma (R/T, S)$ of sections relative to $T$.
This is compatible with realizations: we have
$$
\Gamma (R/T, \kappa (S, n))= \kappa (f_{\ast}(S),n).
$$
Suppose that $f: X\rightarrow Y$ is a morphism of $0$-stacks. Then for a
complex of rational presheaves $C^{\cdot}$ on $X$ the direct image construction
in terms of spectra is the same as the usual higher direct image of complexes
of sheaves (applied to the sheafification of the complex):
$$
f_{\ast}(\sigma (C^{\cdot}/X))= \sigma (({\bf R}f_{\ast}C^{\cdot})/Y).
$$

We extend this just a little bit, in a special case in which it still makes
sense to talk about complexes.  Suppose $X$ is a $1$-stack and $Y$ is a
$0$-stack, with $f: X\rightarrow Y$ a morphism. Suppose $V$ is a local system
of presheaves on $X$ (i.e. for each $Z\in \Xx$, $V(Z)$ is a local system of
rational vector spaces on $X(Z)$). Another way to put this is that $V$ is an
abelian group object over $Z$. We can think of $V$ as being a complex of
presheaves over $X$ (even though we have not defined this notion in general) and
we obtain the spectrum which we denote by $\sigma (V/X)$ over $X$ (even
though this doesn't quite fit in with the general definition of $\sigma$
above), and its realization $\kappa (\sigma (V/X)/X, n)\rightarrow X$ which is
what we would otherwise denote as $K(V/X,n)$.  The higher
direct image ${\bf R}f_{\ast}(V)$ makes sense as a complex of presheaves on $Y$,
and we have the compatibilities
$$
f_{\ast} \sigma (V/X) = \sigma ({\bf R}f_{\ast}(V))
$$
and
$$
\Gamma (X/Y, \kappa (\sigma (V/X)/X, n))= \kappa (\sigma ({\bf
R}f_{\ast}(V)),n).
$$

\begin{proposition}
\label{ComplexOfVB}
Suppose $R$ is an $n$-stack and $S$ is a spectrum over $R$ such that for every
map $Y\rightarrow R$ from a scheme, there is (locally over $Y$ in the etale
topology) a complex of vector bundles $E^{\cdot}_Y$ over $Y$ with $S\times
_RY\cong \sigma (E^{\cdot}_Y/Y)$. Then the realization $\kappa (S/R,n)$ is
geometric over $R$. In particular if $R$ is geoemtric then so is $\kappa
(S/R,n)$. \end{proposition}
{\em Proof:}
In order to prove that the morphism $\kappa (S/R,n)\rightarrow R$ is geometric,
it suffices to prove that for every base change to a scheme $Y\rightarrow R$,
the fiber product $\kappa (S/R, n)\times _RY$ is geometric. But
$$
\kappa (S/R,n)\times _RY= \kappa (\sigma (E^{\cdot}_Y/Y)/Y, n),
$$
so it suffices to prove that for a scheme $Y$ and a complex of vector bundles
$E^{\cdot}$ on $Y$, we have $\kappa (\sigma (E^{\cdot}/Y)/Y, n)$ geometric.

Note that $\kappa (\sigma (E^{\cdot}), n)$ only depends on the part of the
complex
$$
E^0\rightarrow E^1 \rightarrow \ldots \rightarrow E^n\rightarrow E^{n+1}
$$
so we assume that it stops there or earlier. Now we proceed by induction on the
length of the complex.
Define a complex  $F^i= E^{i-1}$ for $i\geq 1$, which has length strictly
smaller than that of $E^{\cdot}$. Let $E^0$ denote the first vector bundle of
$E^{\cdot}$ considered as a complex in degree $0$ only. We have a morphism of
complexes $E^0\rightarrow F^{\cdot}$ and $E^{\cdot}$ is the mapping cone. Thus
$$
\sigma (E^{\cdot}/Y) = \sigma (E^0/Y)\times _{\sigma (F^{\cdot}/Y)}Y
$$
with $Y\rightarrow \sigma (F^{\cdot}/Y)$ the basepoint section. We get
$$
\kappa (\sigma (E^{\cdot}/Y)/Y,n) = K(E^0/Y,n)\times _{\kappa (\sigma
(F^{\cdot}/Y)/Y, n)}Y.
$$
By our induction hypothesis, $\kappa (\sigma
(F^{\cdot}/Y)/Y, n)$ is geometric.  Note that $E_0$ is a smooth group
scheme over
$Y$ so by Lemma \ref{eilenbergExample}, $K(E^0/Y,n)$ is geometric $Y$. By
\ref{fiberprod}, \linebreak  $\kappa (\sigma (E^{\cdot}/Y)/Y,n)$ is geometric.
\eop

{\em Remark:} This proposition is a generalisation to $n$-stacks of
(\cite{LaumonMB} Construction 9.19, Proposition 9.20).  Note that
if $E^{\cdot}$ is a complex where $E^i$ are vector bundles for $i<n$ and $E^n$
is a vector scheme (i.e. something of the form ${\bf V}(\Mm )$ for a coherent
sheaf $\Mm$ in the notation of \cite{LaumonMB}) then we can express $E^n$ as
the kernel of a morphism $U^n\rightarrow U^{n+1}$ of vector bundles (this would
be dual to the presentation of $\Mm$ if we write $E^n = {\bf V}(\Mm )$).
Setting $U^i= E^i$ for $i<n$ we get $\kappa (\sigma (E^{\cdot}), n)=
\kappa (\sigma (U^{\cdot}),n)$. In this way we recover Laumon's and
Moret-Bailly's construction in the case $n=1$.

\begin{corollary}
Suppose $f:X\rightarrow Y$ is a projective flat morphism of schemes, and
suppose that $V$ is a vector bundle on $X$. Then
$\Gamma (X/Y, K(V/X,n))$ is a geometric $n$-stack lying over $Y$.
\end{corollary}
{\em Proof:}
By the discussion at the start of this subsection,
$$
\Gamma (X/Y, K(V/X,n)) = \kappa (\sigma ({\bf R} f_{\ast}(V)/Y)/Y, n).
$$
But by Mumford's method \cite{Mumford}, ${\bf R} f_{\ast}(V)$ is
quasiisomorphic (locally over $Y$) to a complex of vector bundles. By
Proposition \ref{ComplexOfVB} we get that $\Gamma (X/Y, K(V,n))$
is geometric over $Y$.
\eop

Recall that a {\em formal groupoid} is a stack $X_{\Lambda}$ associated to a
groupoid of formal schemes where the object object is a scheme $X$ and the
morphism object is a formal scheme $\Lambda \rightarrow X\times X$ with support
along the diagonal.  We say it is {\em smooth} if the projections
$\Lambda \rightarrow X$ are formally smooth. In this case the cohomology
of the stack $X_{\Lambda}$ with coefficients in vector bundles over
$X_{\Lambda}$ (i.e. vector bundles on $X$ with $\Lambda$-structure meaning
isomorphisms between the two pullbacks to $\Lambda$ satisfying the cocycle
condition on $\Lambda \times _X\Lambda$) is calculated by the {\em de Rham
complex} $\Omega ^{\cdot}_{\Lambda}\otimes _{\Oo}V$ of locally free sheaves
associated to the formal scheme \cite{Illusie} \cite{Berthelot}.

We say that $X_{\Lambda}\rightarrow Y$ is a smooth formal groupoid over $Y$ if
$X_{\Lambda}$ is a smooth fomal groupoid mapping to $Y$ and if $X$ is flat over
$Y$.

\begin{corollary}
Suppose $f: X_{\Lambda} \rightarrow Y$ is a projective smooth formal groupoid
over a scheme $Y$. Suppose that $V$ is a vector bundle on $X_{\Lambda}$ (i.e. a
vector bundle on $X$ with $\Lambda$-structure). Then
$\Gamma (X_{\Lambda}/Y, K(V/X_{\Lambda},n))$ is a geometric $n$-stack lying over
$Y$.
\end{corollary}
{\em Proof:}
By the ``slight extension'' in the discussion at the start of this subsection,
$$
\Gamma (X_{\Lambda}/Y, K(V/X_{\Lambda},n)) = \kappa (\sigma ({\bf R}
f_{\ast}(V)/Y)/Y, n).
$$
But
$$
{\bf R}f_{\ast}(V) = {\bf R} f'_{\ast}(\Omega ^{\cdot}_{\Lambda}\otimes _{\Oo}V)
$$
where $f': X\rightarrow Y$ is the morphism on underlying schemes.
Again  by Mumford's method \cite{Mumford},
${\bf R} f'_{\ast}(\Omega ^{\cdot}_{\Lambda}\otimes _{\Oo}V)$ is
quasiisomorphic (locally over $Y$) to a complex of vector bundles. By the
Proposition \ref{ComplexOfVB} we get that $\Gamma (X_{\Lambda}/Y,
K(V/X_{\Lambda},n))$ is geometric over $Y$.
\eop

\numero{Maps into geometric $n$-stacks}

\begin{theorem}
\label{maps}
Suppose $X\rightarrow S$ is a projective flat morphism. Suppose $T$ is a
connected $n$-stack  which is very presentable (i.e. the fundamental group is
represented by an affine group scheme of finite type denoted $G$ and the higher
homotopy groups are represented by finite dimensional vector spaces).  Then the
morphism $Hom (X/S,T) \rightarrow Bun _G(X/S)= Hom (X/S, BG)$ is a geometric
morphism. In particular $Hom (X/S, T)$ is a locally geometric $n$-stack.
\end{theorem}
{\em Proof:}
Suppose $V$ is a finite dimensional vector space. Let
$$
\Bb (V,n)= BAut (K(V,n))
$$
be the classifying $n+1$-stack for fibrations with fiber $K(V,n)$. It is
connected with fundamental group $GL(V)$ and homotopy group $V$ in dimension
$n+1$ and zero elsewhere.  The truncation morphism
$$
\Bb (V,n)\rightarrow B\, GL(V)
$$
has fiber $K(V,n+1)$ and
admits a canonical section $o: BGL(V)\rightarrow \Bb (V,n)$ (which corresponds
to the trivial fibration with given action of $GL(V)$ on $V$---this fibration
may itself be constructed as $\Bb (V, n-1)$ or in case $n= 2$ as $B(GL(V)
\semidirect V)$).  The fiber of the morphism $o$ is $K(V, n)$, and $BGL(V)$ is
the universal object over $\Bb (V, n)$.

Note that $BGL(V)$ is an geometric $1$-stack (i.e. algebraic stack) and by
Proposition \ref{fibration} applied to the truncation fibration, $\Bb (V, n)$
is a geoemtric $n+1$-stack.

If $X\rightarrow S$ is a projective flat morphism then $Hom (X/S, BGL(V))$ is
a locally geometric $1$-stack (via the theory of Hilbert schemes).   We show
that $p:Hom (X/S, \Bb (V, n))\rightarrow Hom (X/S, BGL(V))$ is a geometric
morphism. For this it suffices to consider a morphism $\zeta :Y\rightarrow Hom
(X/S, BGL(V))$ from a scheme $Y/S$ which in turn corresponds
to a vector bundle $V_{\zeta}$ on $X\times _SY$.  The fiber
of the map $p$ over $\eta$ is $\Gamma (X\times _SY/Y; K(V_{\zeta}, n+1))$
which as we have seen above is geometric over $Y$.  This shows that $p$ is
geometric.  In particular $Hom (X/S, \Bb (V, n))$ is locally geometric.

We now turn to the situation of a general connected geometric and very
presentable $n$-stack $T$.  Consider the truncation morphism $a:T\rightarrow
T':=\tau
_{\leq n-1}T$. We may assume that the theorem is known for  the $n-1$-stack
$T'$. The morphism $a$ is a fibration with fiber $K(V, n)$ so it comes from a
map $b:T' \rightarrow \Bb (V, n)$ and more precisely we have
$$
T = T' \times _{\Bb (V,n)} BGL(V).
$$
Thus
$$
Hom (X/S, T)= Hom (X/S, T') \times _{Hom (X/S,\Bb (V,n))} Hom (X/S,BGL(V)).
$$
But we have just checked that $Hom (X/S,BGL(V))$
and $Hom (X/S,\Bb (V,n))$ are locally geometric, and by hypothesis
$Hom (X/S, T')$ is locally geometric. Therefore by the version of
\ref{fiberprod} for locally geometric $n$-stacks, the fiber product is locally
geometric. This completes the proof. \eop

\begin{theorem}
\label{smoothformal}
Suppose  $(X,\Lambda )\rightarrow S$ is a smooth
projective morphism with smooth formal category structure relative to $S$.
Let $X_{\Lambda}\rightarrow S$ be the resulting family of stacks.
Suppose $T$ is a connected very presentable $n$-stack  which is very presentable
(with fundamental group scheme denoted $G$).  Then the morphism $Hom
(X_{\Lambda }/S,T)  \rightarrow
Hom (X_{\Lambda }/S, BG)$ is a geometric morphism.
In particular $Hom (X_{\Lambda }/S, T)$ is a locally
geometric $n$-stack.
\end{theorem}
{\em Proof:}
The same as before. Note here also that $Hom (X_{\Lambda }/S, BG)$ is an
algebraic stack locally of finite type.
\eop

{\em Remark:} In the above theorems the base $S$ can be assumed to be any
$n$-stack, one looks at morphisms with the required properties when base
changed to any scheme $Y\rightarrow S$.

\subnumero{Semistability}

Suppose $X\rightarrow S$ is a projective flat morphism, with fixed ample class,
and suppose $G$ is an affine algebraic group. We get a notion of semistability
for $G$-bundles (for example, fix the convention that we speak of Gieseker
semistability). Fix also a collection of Chern classes which we denote $c$. We
get a Zariski open substack
$$
Hom ^{\rm se}_c(X/S, BG)\subset Hom (X/S, BG)
$$
(just the moduli $1$-stack of semistable $G$-bundles with Chern classes $c$).
The boundedness property for semistable $G$-bundles with fixed Chern classes
shows that $Hom ^{\rm se}_c(X/S, BG)$ is a geometric $1$-stack.

Now if $T$ is a connected very presentable $n$-stack, let $G$ be the
fundamental group scheme and let $c$ be a choice of Chern classes for
$G$-bundles. Define
$$
Hom ^{\rm se}_c(X/S, T):= Hom (X/S, T)\times _{Hom (X/S, BG)} Hom ^{\rm
se}_c(X/S, BG).
$$
Again it is a Zariski open substack of $Hom (X/S, T)$ and it is a geometric
$n$-stack rather than just locally geometric.

We can do the same in the case of a smooth formal category $X_{\Lambda}
\rightarrow S$.  Make the convention in this case that we ask the Chern classes
to be zero (there is no mathematical need to do this, it is just to conserve
indices, since practically speaking this is the only case we are interested in
below). We obtain a Zariski open substack
$$
Hom ^{\rm se}(X_{\Lambda}/S, BG)\subset Hom (X_{\Lambda}/S, BG),
$$
the moduli stack for semistable $G$-bundles on $X_{\Lambda}$ with vanishing
Chern classes. See \cite{Moduli} for the construction (again the methods given
there suffice for the construction, although stacks are not explicitly
mentionned). Again for any connected very presentable $T$ with fundamental
group scheme $G$ we put
$$
Hom ^{\rm se}(X_{\Lambda}/S, T):= Hom (X_{\Lambda}/S, T)\times _{Hom (X
_{\Lambda}/S, BG)} Hom ^{\rm
se}(X_{\Lambda}/S, BG).
$$
It is a geometric $n$-stack.

Finally we note that in the case of the relative de Rham formal category
$X_{DR/S}$ semistability of principal $G$-bundles is automatic (as is the
vanishing of the Chern classes). Thus
$$
Hom ^{\rm se}(X_{DR/S}/S, T)= Hom (X_{DR/S}/S, T)
$$
and $Hom (X_{DR/S}, T)$ is already a geometric $n$-stack.

\subnumero{The Brill-Noether locus}

Suppose $G$ is an algebraic group and $V$ is a representation.  Define the
$n$-stack $\kappa (G,V,n)$ as the fibration over $K(G,1)$ with fiber
$K(V,n)$ where $G$ acts on $V$ by the given representation and such that
there is a section.  Let $X$ be a projective variety. We have a morphism
$$
Hom (X, \kappa (G,V,n))\rightarrow Hom (X, K(G,1))= Bun _G(X).
$$
The fiber over a point $S\rightarrow Bun _G(X)$ corresponding to a principal
$G$-bundle $P$ on $X\times S$ is the relative section space
$$
\Gamma (X\times S/S, K(P\times ^GV/X\times S, n)).
$$
By the compatibilities given at the start of the section on Dold-Puppe, this
relative section space is the $n$-stack corresponding to the  direct image
$Rp_{1,\ast}(P\times ^GV)$  which is a complex over $S$.   Note that this
complex is quasiisomorphic to a complex of vector bundles.  Thus we have:

\begin{corollary}
\label{BN}
The morphism
$$
Hom (X, \kappa (G,V,n))\rightarrow  Bun _G(X)
$$
is a morphism of geometric $n$-stacks.
\end{corollary}
\eop

{\em Remark:}  The $Spec (\cc )$-valued points of $Hom (X, \kappa
(G,V,n))$ are the pairs $(P, \eta )$ where $P$ is a principal $G$-bundle
on $X$ and $\eta \in H^n(X, P\times ^GV)$.

Thus $Hom (X, \kappa
(G,V,n))$ is a geometric $n$-stack whose $Spec (\cc )$-points are the
Brill-Noether set of vector bundles with cohomology classes on $X$.

\subnumero{Some conjectures}

We give here some conjectures about the possible extension of the above results
to any (not necessarily connected) geometric $n$-stacks $T$.

\begin{conjecture}
If $T$ is a geometric $n$-stack which is very presentable in the sense of
\cite{RelativeLie} (i.e. the fundamental groups over artinian base are affine,
and the higher homotopy groups are vector sheaves) then for any smooth (or just
flat?) projective morphism $X\rightarrow S$ we have that $Hom (X/S, T)$ is
locally geometric. \end{conjecture}

\begin{conjecture}
\label{KGm2}
If $T= K({\bf G}_m , 2)$ then for a flat projective morphism $X\rightarrow S$,
$Hom (X/S, T)$ is locally geometric. Similarly if $G$ is {\em any} group scheme
of finite type (e.g. an abelian variety)
then $Hom (X/S, BG)$ is locally geometric.
\end{conjecture}

Putting together with the previous conjecture we can make:

\begin{conjecture}
If $T$ is a geoemtric $n$-stack whose $\pi _i$ are vector sheaves for $i\geq 3$
then $Hom (X/S, T)$ is locally geometric.
\end{conjecture}

Note that Conjecture \ref{KGm2} cannot be true if $K({\bf G}_m, 2)$ is
replaced by $K({\bf G}_m, i)$ for $i\geq 3$, for in that case the morphism
stacks will themselves be only locally of finite type. Instead we will get a
``slightly geometric'' $n$-stack as discussed in \S 3.  One could make the
following conjecture:

\begin{conjecture}
If $T$ is any geometric (or even locally or slightly geometric) $n$-stack and
$X\rightarrow S$ is a flat projective morphism then $Hom (X/S, T)$ is slightly
geometric.
\end{conjecture}

After these somewhat improbable-sounding conjectures, let finish by making
a more
reasonable statement:
\begin{conjecture}
If $T$ is a very presentable geometric $n$-stack and $X$ is a smooth projective
variety then $Hom (X_{DR}, T)$ is again geometric.
\end{conjecture}

Here, we have already announced the finite-type result in the statement that
\linebreak
$Hom (X_{DR}, T)$ is very presentable \cite{kobe} (I have not yet circulated
the proof, still checking the details...).

\subnumero{GAGA}

Let $\Xx ^{\rm an}$ be the site of complex analytic spaces with the etale (or
usual--its the same) topology.  We can make similar definitions of geometric
$n$-stack on $\Xx ^{\rm an}$ which we will now denote by {\em analytic
$n$-stack} (in case of confusion...).  There are similar definitions of
smoothness and so on.

There is a morphism of sites from the analytic to the algebraic sites.

If $T$ is a geometric $n$-stack on $\Xx$ then its pullback by this morphism (cf
\cite{realization}) is an analytic $n$-stack which we denote by $T^{\rm an}$.

We have:

\begin{theorem}
\label{gaga}
Suppose $T$ is a connected very presentable geometric $n$-stack.  Suppose
$X\rightarrow S$ is a flat projective morphism (resp. suppose
$X_{\Lambda}\rightarrow S$ is the morphism associated to a smooth formal
category over $S$). Then the natural morphism
$$
Hom (X/S, T)^{\rm an} \rightarrow Hom (X^{\rm an}/S^{\rm an}, T^{\rm an})
$$
$$
\left( \mbox{resp.} Hom (X_{\Lambda}/S, T)^{\rm an} \rightarrow Hom
(X_{\Lambda}^{\rm
an}/S^{\rm an}, T^{\rm an})
\right)
$$
is an isomorphism of analytic $n$-stacks.
\end{theorem}
{\em Proof:}
Just following through the proof of the facts that $Hom (X/S, T)$
or $Hom (X_{\Lambda}/S, T)$ are geometric, we can keep track of the analytic
case too and see that the morphisms are isomorphisms along with the main
induction.
\eop

{\em Remarks:}
\newline
(1)\, This GAGA theorem holds for $X_{DR}$ with coefficients in any very
presentable $T$ (not necessarily connected) \cite{kobe}.
\newline
(2)\, In \cite{kobe} we also give a ``GFGA'' theorem for $X_{DR}$ with
coefficients in a  very presentable $n$-stack.
\newline
(3)\, The GAGA theorem does not hold with coefficients in $T= K({\bf G}_m ,
2)$. Thus the condition that the higher homotopy group sheaves of $T$ be vector
sheaves is essential. Maybe it could be weakened by requiring just that the
fibers over artinian base schemes be unipotent (but this might also be
equivalent to the vector sheaf condition). \newline
(4)\, Similarly the GAGA theorem does not hold with coefficients in
$T= BA$ for an abelian variety $A$; thus again the hypothesis that the fibers of
the fundamental group sheaf over artinian base be affine group schemes, is
essential.

\numero{The tangent spectrum}

We can treat a fairly simple case of the conjectures outlined above: maps from
the spectrum of an Artin local algebra of finite type.

\begin{theorem}
\label{mapsFromArtinian}
Let $X=Spec (A)$ where
$A$ is artinian, local, and of finite type over $k$. Suppose $T$ is a
geometric $n$-stack. Then $Hom (X,T)$ is a geometric $n$-stack. If
$T\rightarrow T'$ is a geometric smooth morphism of $n$-stacks then $Hom (X,
T)\rightarrow Hom (X, T')$ is a smooth geometric morphism of $n$-stacks.
\end{theorem}
{\em Proof:}
We prove the following statement: if $Y$ is a scheme and $A$ as in the theorem,
and if $T\rightarrow Y\times Spec (A)$ is a geometric (resp. smooth geometric)
morphism of $n$-stacks then  $\Gamma (Y\times Spec (A)/Y, T)$ is geometric
(resp.
smooth geometric) over $Y$. The proof is by induction on $n$; note that it
works for $n=0$. Now in general choose a smooth surjection $X\rightarrow T$
from a scheme.  Then $\Gamma (Y\times Spec (A)/Y, X)$
is a scheme over $Y$, and if $X$ is smooth over $Y$ then the section scheme is
smooth over $Y$. We have a surjection
$$
a:\Gamma (Y\times Spec (A)/Y, X)\rightarrow
\Gamma (Y\times Spec (A)/Y, T),
$$
and for $Z\rightarrow \Gamma (Y\times Spec (A)/Y, X)$ (which amounts to
a section morphism $Z\times Spec (A)\rightarrow T$) the fiber product

$$
\Gamma (Y\times Spec (A)/Y, X)\times _{ \Gamma (Y\times Spec (A)/Y, T)}
Z
$$
is equal to
$$
\Gamma (Z\times Spec (A)/Z, X\times _T(Z\times Spec (A))).
$$
But $X\times _T(Z\times Spec (A)$ is a smooth $n-1$-stack over
$Z\times Spec (A)$ so by induction this section stack is geometric and smooth
over $Z$.  Thus our surjection $a$ is a smooth geometric morphism so
$\Gamma (Y\times Spec (A)/Y, T)$ is geometric.  The smoothness statement
follows immediately.
\eop

We apply this to define the {\em tangent spectrum} of a  geometric
$n$-stack. This is a generalization of the discussion at the end of
(\cite{LaumonMB} \S 9), although we use a different approach because I
don't have
the courage to talk about cotangent complexes!

Recall from \cite{Adams} Segal's infinite loop space machine: let $\Gamma$ be
the category whose objects are finite sets and where the morphisms from
$\sigma$ to $\tau$ are maps $P(\sigma )\rightarrow P(\tau )$ preserving
disjoint unions (here $P(\sigma )$ is the set of subsets of $\sigma$).
A morphism is determined, in fact, by the map $\sigma \rightarrow P(\tau )$
taking different elements of $\sigma$ to disjoint subsets of $\tau$ (note that
the empty set must go to the empty set).  Let $[n]$
denote the set with $n$ elements. There is a functor $s:\Delta \rightarrow
\Gamma$ sending the the ordered set $\{ 0,\ldots , n\}$ to the finite set $\{ 1,
\ldots , n\}$---see \cite{Adams} p. 64 for the formulas for the morphisms.

Segal's version of an
$E_{\infty}$-space (i.e. infinite loop space) is a contravariant functor
$\Psi : \Gamma \rightarrow Top$ such that the associated  simplicial
space   (the composition $\Psi \circ s$) satisfies Segal's condition
\cite{Adams}. In order to really get an infinite loop space it is also required
that $\Psi (\emptyset )$ be a point (although this condition seems to have been
lost in Adams' very brief treatment).

Segal's machine is then a classifying space functor $B$ from
special $\Gamma$-spaces to special $\Gamma$-spaces. This actually works even
without the condition that $\Phi (\emptyset )$ be a point, however the
classifying space construction is the inverse to the {\em relative} loop space
construction over $\Phi (\emptyset )$.  Note that since $\emptyset$ is a final
object in $\Gamma$ the components of a $\Gamma$-space are provided with a
section from $\Phi (\emptyset )$. If
$\Phi$ is a special $\Gamma$-space then $B^n\Phi$ is again a special $\Gamma$
space with
$$
B^n\Phi (\emptyset )= \Phi (\emptyset )
$$
and
$$
\Omega ^n(B^n\Phi ([1])/\Phi (\emptyset ))= \Phi ([1]).
$$
The notion of $\Gamma$-space (say with $\Phi [1] $ rational over $\Phi
(\emptyset )=R$) is another replacement for our notion of spectrum over $R$; we
get to our notion as defined above by looking at $B^N\Phi ([1])$.

The above discussion makes sense in the context of presheaves of spaces over
$\Xx$ hence in the context of $n$-stacks.

We now try to apply this in our situation to construct the tangent spectrum.
For any object $\sigma \in \Gamma$ let ${\bf A}^\sigma$ be the affine space
over $k$
with basis the set $\sigma$.  An element of ${\bf A}^\sigma$ can be written as
$\sum _{i\in \sigma} a_ie_i$ where $e_i$ are the basis elements and $a_i\in k$.
Given a map $f:\sigma \rightarrow P(\tau )$ we define a map
$$
{\bf A}^f : {\bf A}^{\sigma} \rightarrow {\bf A}^{\tau}
$$
$$
\sum a_i e_i \mapsto \sum _{i\in \sigma} \sum _{j\in f(i)\subset \tau}
a_i e_j.
$$
For example there are four morphisms  from $[1]$ to $[2]$, sending $1$ to
$\emptyset$, $\{ 1\}$, $\{ 2\}$ and $\{ 1,2\}$ respectively.  These correspond
to the constant morphism, the two coordinate axes, and the diagonal from ${\bf
A}^1$ to ${\bf A}^2$.  We get a covariant functor from $\Gamma$ to the category
of affine schemes.

For a finite set $\sigma$ let $D^{\sigma}$ denote the subscheme of ${\bf
A}^{\sigma}$ defined by the square of the maximal ideal defining the origin.
These fit together into a covariant functor from $\Gamma$ to the category of
artinian local schemes of finite type over $k$.

If $T$ is a geometric  $n$-stack thought of as a strict presheaf of spaces, then
the functor
$$
\Theta :\sigma \mapsto Hom (D^{\sigma}, T)
$$
is a contravariant functor from $\Gamma$ to the category of
geometric $n$-stacks, with $\Theta (\emptyset )=T$. To see that it satisfies
Segal's condition we have to check that the map
$$
Hom (D^n, T)\rightarrow Hom (D^1, R)\times _R \ldots \times _THom (D^1, T)
$$
is an equivalence.  Once this is checked we obtain a spectrum over $T$ whose
interpretation in our terms is as the $N$-stack $B^N\Phi ([1])$.

In the statement of the following theorem we will normalize our relationship
between complexes and spectra in a different way from before---the most natural
way for our present purposes.

\begin{theorem}
\label{tangent}
Suppose $T$ is a geometric $n$-stack. The above construction gives a spectrum
$\Theta (T)\rightarrow T$ which we call the {\em tangent spectrum of $T$}.
If $Y\rightarrow T$ is a morphism from a scheme then $\Theta (T)\times _TY$
is equivalent to $\sigma (E^{\cdot}/Y)$ for a complex
$$
E^{-n}\rightarrow \ldots \rightarrow E^0
$$
with $E^i$ vector bundles ($i<0$) and $E^0$ a vector scheme over $Y$.
Furthermore if $T$ is smooth then $E^0$ can be assumed to be a vector bundle.
In particular, $\kappa (\Theta (T)/T, n)$ is geometric, and if $T$ is smooth
then $\Theta (T)$ is geometric.
\end{theorem}
{\em Proof of \ref{tangent}:}
The first task is to check the above condition for $\Theta$ to be a special
$\Gamma$-space. Suppose in general that $A,B\subset C$ are closed artinian
subschemes of an artinian scheme with the extension property that for any
scheme $Y$ the morphisms from $C$ to $Y$ are the same as the pairs of morphisms
$A,B\rightarrow Y$ agreeing on $A\cap B$.  We would like to show that for any
geometric stack $T$,
$$
Hom (C,T)\rightarrow Hom (A,T)\times _{Hom(A\cap B, T)}Hom (B,T)
$$
is an equivalence.  We have a similar relative statement for sections of a
geometric morphism $T\rightarrow Y\times C$ for a scheme $Y$.  We prove the
relative statement by induction on the truncation level $n$, but for
simplicity use the notation of the absolute statement.  Let $X\rightarrow
T$ be a
smooth geometric morphism from a scheme.  Then consider the diagram
$$
\begin{array}{ccc}
Hom (C,X) &\stackrel{\cong}{\rightarrow}&
Hom (A,X)\times _{Hom(A\cap B, X)}Hom (B,X) \\
\downarrow && \downarrow \\
Hom (C,T)&\rightarrow &Hom (A,T)\times _{Hom(A\cap B, T)}Hom (B,T).
\end{array}
$$
It suffices to prove that for a map from a scheme $Y\rightarrow
Hom (C,T)$ the morphism on fibers is an equivalence. The fiber on the left is
$$
Hom (C,X)\times _{Hom (C,T)}Y= \Gamma (Y\times C, X\times _{T}(Y\times C)),
$$
whereas the fiber on the right is
$$
\Gamma (Y\times A, X\times _{T}(Y\times A))
\times _{\Gamma (Y\times (A\cap B), X\times _{T}(Y\times (A\cap B)))}
\Gamma (Y\times B, X\times _{T}(Y\times B)).
$$
By the relative version of the statement for the $n-1$-stack
$X\times _{T}(Y\times C)$ over $Y\times C$, the map of fibers is an equivalence,
so the map
$$
Hom (C,T)\rightarrow Hom (A,T)\times _{Hom(A\cap B, T)}Hom (B,T)
$$
is an equivalence.

Apply this inductively with $C= D^n$, $A= D^1$ and $B= D^{n-1}$ (so $A\cap
B=D^0$).  We obtain the required statement, showing that $\Theta$ is a special
$\Gamma$-space relative to $T$.  It integrates to a spectrum which we denote
$\Theta (T)\rightarrow T$.

Note that if $T=X$ is a scheme considered as an $n$-stack then $\Theta (X)$ is
just the spectrum associated to the complex consisting of the tangent vector
scheme of $X$ in degree $0$.  We obtain the desired statement in this case.

If $R\rightarrow T$ is a morphism of geometric $n$-stacks then
we obtain a morphism of spectra
$$
\Theta (R) \rightarrow  \Theta (T)\times _T R .
$$
The cofiber (i.e. $B$ of the fiber) we denote by $\Theta (R/T)$.
We prove more generally---by induction on $n$---that if $T\rightarrow Y$ is a
geometric morphism from an $n$-stack to a scheme, and if $Y\rightarrow T$ is a
section then $\Theta (T/Y)\times _TY$ is associated (locally on $Y$) to a
complex
of vector bundles and a vector scheme at the end; with the last vector scheme
being a bundle if the morphism is smooth.  Note that it is true for $n=0$.  For
any $n$ choose a smooth geometric morphism $X\rightarrow T$ and we may assume
(by etale localization) that there is a lifting of the section to $Y\rightarrow
X$. Now there is a triangle of spectra (i.e. associated to a triangle of
complexes in the derived category)
$$
\Theta (X)\times _XY \rightarrow \Theta (T)\times _TY \rightarrow B\Theta
(X/T)\times _XY.
$$
On the other hand,
$$
B\Theta (X/T)\times _XY=B\Theta (X\times _TY/Y)\times _{X\times _TY}Y.
$$
By induction this is associated to a complex as desired, and we know already
that $\Theta (X)\times _XY$ is associated to a complex as desired.  Therefore
$\Theta(T)\times _TY$ is an extension of complexes of the desired form, so it
has the desired form. Note that since $X\times _TY\rightarrow Y$ is smooth,
by the induction hypothesis we get that $B\Theta (X/T)\times _XY$ is associated
to a complex of bundles.

If the morphism $T\rightarrow Y$ is smooth then the last term in the complex
will be a bundle (again following through the same induction).
\eop

If $T$ is a smooth geometric $n$-stack and $P: Spec (k)\rightarrow T$ is a
point then we say that the {\em dimension of $T$ at $P$} is the alternating sum
of the dimensions of the vector spaces in the complex making up the
complex associated to $P^{\ast} (\Theta (T))$. This could, of course, be
negative.

For example if $G$ is an algebraic group then the dimension of $BG$ at any
point is $-dim (G)$.
More generally if $A$ is an abelian group scheme smooth over a base $Y$ then
$$
dim (K(A/Y, n))= dim (Y) + (-1)^ndim (A).
$$

\numero{De Rham theory}

We will use certain geometric $n$-stacks as coefficients to look at the de Rham
theory of a smooth projective variety.  The answers come out to be  geometric
$n$-stacks. (One could also try to look at de Rham theory {\em for} geometric
$n$-stacks, a very interesting problem but not what is meant by the title of
the present section).

If $X$ is a smooth projective variety let $X_{DR}$ be the stack (which is
actually a sheaf of sets) associated to the formal category whose object object
is $X$ and whose morphism object is the formal completion of the diagonal in
$X\times X$.  Another cheaper definition is just to say
$$
X_{DR}(Y):= X(Y^{\rm red}).
$$
If $f:X\rightarrow S$ is a smooth morphism, let
$$
X_{DR/S}:= X_{DR}\times _{S_{DR}}S.
$$
It is the stack associated to a smooth formal groupoid over $S$ (take the formal
completion of the diagonal in $X\times _SX$).

The cohomology of $X_{DR}$ with coefficients in an algebraic group scheme is
the same as the de Rham cohomology of $X$ with those coefficients.
We treat this in the case of coefficients in a vector space, or in case of
$H^1$ with coefficients in an affine group scheme. Actually the statement is a
more general one about formal categories. Suppose $(X,\Lambda )\rightarrow S$
is a smooth formal groupoid over $S$  which
we can think of as a smooth scheme $X/S$ with a formal scheme $\Lambda$ mapping
to $X\times _SX$ and provided with an associative product structure. There
is an associated {\em de Rham complex} $\Omega ^{\cdot} _{\Lambda}$ on $X$
(cf \cite{Berthelot} \cite{Illusie})---whose components are locally free
sheaves on $X$ and where the differentials are first order differential
operators.  Let $X_{\Lambda}$ denote the stack associated to the formal
groupoid. It is the stack associated to the presheaf of groupoids which to $Y\in
\Xx$ associates the groupoid whose objects are $X(Y)$ and whose morphisms are
$\Lambda (Y)$.

Suppose $V$ is a vector bundle over $X_{\Lambda}$, that is a vector bundle
on $X$ together with isomorphisms $p_1^{\ast} V\cong p_2^{\ast} V$ on $\Lambda$
satisfying the cocycle condition on $\Lambda \times _X \Lambda$.
We can define the cohomology sheaves on $S$, $H^i(X_{\Lambda}/S, V)$ which will
be equal to $\pi _0(\Gamma (X_{\Lambda }/S; K(V, i))$ in our notations. These
cohomology sheaves can be calculated using the de Rham complex: there is a
twisted de Rham complex $\Omega ^{\cdot}_{\Lambda} \otimes _{\Oo}V$ whose
hypercohomology is $H^i(X_{\Lambda}/S, V)$.

When applied to the de Rham formal category (the trivial example introduced in
\cite{Berthelot} in characteristic zero) whose associated stack is the sheaf of
sets $X_{DR/S}$, we obtain the usual de Rham complex $\Omega ^{\cdot}_{X/S}$
relative to $S$.  A vector bundle $V$ over $X_{DR/S}$ is the same thing as a
vector bundle on $X$ with integrable connection, and the twisted de Rham
complex is the usual one.  Thus in this case we have
$$
\pi _0(\Gamma (X_{DR/S}/S, K(V,i)))= {\bf H}^i(X/S, \Omega ^{\cdot}_{X/S}\otimes
V).
$$
We can describe more precisely  $\Gamma (X_{DR/S}/S, K(V,i))$
as being the$i$-stack obtained by applying Dold-Puppe to the right derived
direct
image complex $Rf _{\ast} (\Omega ^{\cdot}_{X/S}\otimes
V)[i]$ (appropriately shifted).

For the first cohomology with coefficients in an affine algebraic group $G$,
note that a principal $G$-bundle on $X_{DR}$ is the same thing as a principal
$G$-bundle on $X$ with integrable connection. We have that the $1$-stack
$\Gamma (X_{DR/S}/S, BG)$ on $S$ is the moduli stack of principal $G$-bundles
with relative integrable connection on $X$ over $S$. For $X\rightarrow S$
projective this is constructed in \cite{Moduli} (in fact, there we construct
the representation scheme of framed principal bundles; the moduli stack is
immediately obtained as an algebraic stack, the quotient stack by the action
of $G$ on the scheme of framed bundles).

Of course we have seen in \ref{smoothformal} that for any smooth formal category
$(X,\Lambda )$ over $S$ and any connected very presentable $n$-stack $T$,
the morphism $n$-stack $Hom (X_{\Lambda} /S, T)$ is a locally geometric
$n$-stack.   Recall that we have defined the {\em semistable} morphism stack
$Hom ^{\rm se}(X_{\Lambda} /S, T)$ which is geometric; but in our case all
morphisms $X_{DR/S}\rightarrow BG$ (i.e. all principal $G$-bundles with
integrable connection) are semistable, so in this case we find that $Hom
(X_{DR/S}/S, T)$ is a geometric $n$-stack. In fact it is just a successive
combination of the above discussions applied according to the Postnikov
decomposition of $T$.

\subnumero{De Rham theory on the analytic site}
The same construction works for smooth objects in the analytic site.
Suppose $f:X\rightarrow S$ is a smooth analytic morphism.
Here
we would like to consider any connected $n$-stack $R$ whose homotopy
groups are represented by analytic Lie groups.  Such an $R$ is automatically an
analytic $n$-stack (by the analytic analogue of \ref{fibration}).  We call these
the ``good connected analytic $n$-stacks'' since we haven't yet proven that
every
connected analytic $n$-stack must be of this form (I suspect that to be true but
don't have an argument).

If $G$ is an analytic Lie group, a map $X_{DR/S}\rightarrow
BG$ is a principal
$G$-bundle $P$ on $X$ together with an integrable connection
relative to $S$.

Suppose $A$ is an analytic abelian Lie group with action of $G$.  Then we can
form the analytic $n$-stack  $\kappa (G, A, n)$ with fundamental group $G$ and
$\pi _n= A$.  Given a map $X_{DR/S}\rightarrow BG$ corresponding to a
principal bundle $P$, we would like to study the liftings into $\kappa (G,A,n)$.
We obtain the twisted analytic Lie group $A_P:= P\times ^GA$ over $X$ with
integrable connection relative to $S$.  Let $V$ denote the universal
covering group of $A$ (isomorphic to $Lie (A)$, thus $G$ acts here) and let
$L\subset V$ denote the kernel of the map to $A$. Note that $V$ is a complex
vector space and  $L$ is a lattice isomorphic to $\pi _1(A)$. Again $G$ acts on
$L$. We obtain after twisting $V_P$ and $L_P$.  Note that $V_P$ is provided
with an integrable connection relative to $S$.  The following Deligne-type
complex calculates the cohomology of $A_P$:
$$
C^{\cdot}_{\Dd} (A_P):= \{ L_P \rightarrow V_P \rightarrow \Omega
^1_{X/S}\otimes
_{\Oo} V_P \rightarrow \ldots \} .
$$
The $n$-stack $\Gamma (X_{DR/S}/S, K(A_P,n))$ is again obtained by applying
Dold-Puppe to the shifted  right derived direct image complex
$Rf_{\ast}(C^{\cdot}_{\Dd}(A_P))[n]$.   We can write
$ C^{\cdot}_{\Dd} (A_P)$ as the mapping cone of a map of complexes
$L_P \rightarrow U^{\cdot}_P$.

If $f$ is a projective morphism then applying GAGA and the argument of Mumford
(actually I think there is an argument of Grauert which treats this for any
proper map), we get that
$$
Rf_{\ast}(C^{\cdot}_{\Dd}(U^{\cdot}_P))
$$
is quasiisomorphic to a complex of analytic Lie groups (vector bundles in this
case).  On the other hand, locally on the base the direct image
$Rf_{\ast}(C^{\cdot}_{\Dd}(L_P)$ is a trivial complex so quasiisomorphic to a
complex of (discrete) analytic Lie groups. The direct image
$Rf_{\ast}(C^{\cdot}_{\Dd}(A_P))$ is the mapping cone of a map of these
complexes, so the associated spectrum fits into a fibration sequence. The base
and the fiber are analytic $N$-stack so the total space is also an analytic
$N$-stack. Thus the spectrum associated to
$Rf_{\ast}(C^{\cdot}_{\Dd}(A_P))$ is analytic over $S$. In particular its
realization $\Gamma (X_{DR/S}/S, K(A_P,n))$ is a
geometric $n$-stack over $S$.

For $Hom (X_{DR/S}/S, BG)$ we can use the Riemann-Hilbert correspondence (see
below) to see that it is an analytic $1$-stack.
The same argument as in Theorem \ref{maps} now shows
that for any good connected analytic $n$-stack $T$, the $n$-stack of morphisms
$Hom (X_{DR/S}/S, T)$ is an analytic $n$-stack over $S$.

If the base $S$ is a point we don't need to make use of Mumford's argument,
so the same holds true for any proper smooth analytic space $X$.

{\em Caution:}  There is (at least) one gaping hole in the above argument,
because we are applying Dold-Puppe for complexes of $\zz$-modules such as $L_P$
or its higher direct image, which are not complexes of rational vector spaces.
Thus this doesn't fit into the previous discussion of Dold-Puppe, spectra etc.
as we have set it up.  In particular there may be problems with torsion, finite
groups or subgroups of finite index in the above discussion. The reader is
invited to try to figure out how to fill this in (and to let me know if he
does).

\subnumero{The Riemann-Hilbert correspondence}

We can extend to our cohomology stacks the classical Riemann-Hilbert
correspondence.  We start with a statement purely in the analytic case.
In order to avoid confusion between the analytic situation and the algebraic
one, we will append the superscript {\em `an'} to objects in the analytic site,
even if they don't come from objects in the algebraic site. We will make clear
in the hypothesis whenever our analytic objects actually come from algebraic
ones.

\begin{theorem}
\label{analyticRiemannHilbert}
Suppose $T^{\rm an}$ is a good connected analytic $n$-stack, and suppose
$X^{\rm an}$ is a smooth proper complex analytic space. Define $X^{\rm an}_{DR}$
as above. Let $X^{\rm an}_B$ denote the $n$-stack associated to the constant
presheaf of spaces which to each $Y^{\rm an}$ associates the topological space
$X^{\rm top}$. Then there is a natural equivalence of analytic $n$-stacks $Hom
(X^{\rm an}_{DR}, T^{\rm an}) \cong Hom (X^{\rm an}_B, T^{\rm an})$.
\end{theorem}
{\em Proof:}
By following the same outline as the argument given in \ref{maps}, it
suffices to
see this for the cases $T^{\rm an} = BG^{\rm an}$ for an analytic Lie group
$G^{\rm an}$, and $T^{\rm an}= \Bb (A^{\rm an}, n)$ for an abelian analytic Lie
group  $A^{\rm an}$. In the
second case we reduce to the case of cohomology with  coefficients in a twisted
version of $A^{\rm an}$.  We now leave it to the reader to verify these cases
(which are standard examples of using analytic de Rham cohomology to calculate
singular cohomology).
\eop

{\em Remark:}  For convenience we have stated only the absolute version. We
leave it to the reader to obtain a relative version for a smooth projective
morphism $f: X\rightarrow S$.

Now we turn to the algebraic situation.
We can combine the above result with GAGA  to obtain:

\begin{theorem}
\label{algebraicRiemannHilbert}
Suppose $T$ is a connected very  presentable algebraic $n$-stack, and suppose
$X$ is a smooth projective variety. Define $X_{DR}$ as above.
Let $X_B$ denote the $n$-stack associated to the constant presheaf of spaces
which to each $Y$ associates the topological space $X^{\rm top}$. Then
there is a
natural equivalence of analytic $n$-stacks
$Hom (X_{DR}, T)^{\rm an} \cong Hom (X_B, T)^{\rm an}$.
\end{theorem}
{\em Proof:}
By GAGA we have
$$
Hom (X_{DR}, T)^{\rm an} \cong Hom (X^{\rm an}_{DR}, T^{\rm an}).
$$
Similarly the calculation of $Hom (X_B,T)$ using a cell decomposition of
$X_B$ and fiber products yields the equivalence
$$
Hom (X_B,T)^{\rm an}\cong Hom (X^{\rm an}_B, T^{\rm an}).
$$
Putting these together we obtain the desired equivalence.
\eop

\subnumero{The Hodge filtration}

Let $H:= {\bf A}^1/{\bf G}_m$ be the quotient stack of the affine line modulo
the action of the multiplicative group. This has a Zariski open substack which
we denote $1\subset H$; note that $1\cong Spec (\cc )$.  There is a closed
substack $0\subset H$ with $0\cong B{\bf G}_m$.

As in \cite{SantaCruz} we can define a smooth formal category $X_{\rm
Hod}\rightarrow H$ whose fiber over $1$ is $X_{DR}$ and whose fiber over $0$ is
$X_{Dol}$.

Suppose $T$ is a connected very presentable $n$-stack. Then we obtain the
relative semistable morphism stack
$$
Hom ^{\rm se}(X_{\rm  Hod}/H, T) \rightarrow H.
$$
In the case $T=BG$ this was interpreted as the {\em Hodge filtration on $\Mm
_{DR}=Hom (X_{DR}, BG)$}. Following this interpretation, for any connected very
presentable $T$ we call this relative morphism stack the {\em Hodge filtration
on the higher nonabelian cohomology stack $Hom (X_{DR}, T)$}.

Note that when $T= K(\Oo  , n)$ we recover the usual Hodge filtration on the
algebraic de Rham cohomology, i.e. the cohomology of $X_{DR}$ with coefficients
in $\Oo$.

The above general definition is essentially just a mixture of the case $BG$ and
the cases $K(\Oo , n)$ but possibly with various twistings.

{\em The analytic case:}  The above discussion works equally well for a smooth
proper analytic variety $X$. For any good connected analytic $n$-stack $T$ we
obtain the relative morphism stack
$$
Hom (X_{\rm  Hod}/H^{\rm an}, T) \rightarrow H^{\rm an}.
$$
Note that there is no question of semistability here. The moduli stack of flat
principal $G$-bundles $Hom (X_{\rm Hod}/H^{\rm an}, BG)$ is still an analytic
$n$-stack because in the analytic category there is no distinction between
finite type and locally finite type.

In case $X$ is projective and $G= \pi _1(T)$ affine algebraic we can
put in the semistability condition and get
$$
Hom ^{\rm se}(X_{\rm  Hod}/H^{\rm an}, T) \rightarrow H^{\rm an}.
$$
If $T$ is the analytic stack associated to an algebraic geometric $n$-stack
then this analytic morphism stack is the analytic stack associated to the
algebraic morphism stack.

\subnumero{The Gauss-Manin connection}
Suppose $X\rightarrow S$ is a smooth projective morphism and $T$ a connected
very presentable $n$-stack.  The formal category $X_{DR/S}\rightarrow S$ is
pulled back from the morphism $X_{DR}\rightarrow S_{DR}$ via the map
$S\rightarrow S_{DR}$. Thus
$$
Hom (X_{DR/S}/S,T) = Hom (X_{DR}/S_{DR}, T)\times _{S_{DR}}S.
$$
Thus the morphism stack $Hom (X_{DR/S}/S,T)$ descends down to an $n$-stack over
$S_{DR}$.  If $Y\rightarrow S_{DR}$ is a morphism from a scheme then locally in
the etale topology it lifts to $Y\rightarrow S$.  We have
$$
Hom (X_{DR}/S_{DR},T)\times _{S_{DR}}Y=
Hom (X_{DR}\times _{S_{DR}}Y/Y, T)=
$$
$$
Hom (X_{DR/S}\times _SY/Y,T)= Hom
(X_{DR/S}/S,T)\times _SY.
$$
The right hand side is a geometric $n$-stack, so this shows that the morphism
$$
Hom (X_{DR}/S_{DR},T)\rightarrow S_{DR}
$$
is geometric.  This descended structure is the {\em Gauss-Manin connection} on
$Hom (X_{DR/S}/S,T)$.  In the case $T=BG$ this gives the Gauss-Manin connection
on the moduli stack of $G$-bundles with flat connection (cf \cite{Moduli},
\cite{SantaCruz}). In the case $T= K(V,n)$ this gives the Gauss-Manin
connection on algebraic de Rham cohomology.

In \cite{SantaCruz} we have indicated, for the case $T=BG$, how to obtain the
analogues of {\em Griffiths transversality} and {\em regularity} for the
Hodge filtration and Gauss-Manin connection.
Exactly the same constructions work
here. We briefly review how this works.  Suppose $X\rightarrow S$ is a
smooth projective family
over a quasiprojective base (smooth, let's say) which extends to a family
$\overline{X}\rightarrow \overline{S}$ over a normal crossings compactification
of the base.  Let $D= \overline{X}-X$ and $E= \overline{S}-S$.
Recall that $\overline{X}_{\rm Hod}(\log D)\rightarrow H$ is the
smooth formal category whose underlying space (stack, really, since we have
replaced ${\bf A}^1$ by its quotient $H$) is $X\times H$ and whose associated
de Rham complex is $(\Omega ^{\cdot}_{\overline{X}}(\log D), \lambda d)$ where
$\lambda $ is the coordinate on $H$ (actually to be correct we have to twist
everything by line bundles on $H$ to reflect the quotient by ${\bf G}_m$ but I
won't put this into the notation). Similarly we obtain the formal category
$\overline{S}_{\rm Hod}(\log E)\rightarrow H$, with a morphism
$$
\overline{X}_{\rm Hod}(\log D) \rightarrow
\overline{S}_{\rm Hod}(\log E).
$$
If we pull back by $\overline{S}\rightarrow \overline{S}_{\rm Hod}(\log E)$
then we get a smooth formal category over $\overline{S}$. Thus by
\ref{smoothformal} for any connected very presentable  $n$-stack $T$  the
morphism $$
Hom (\overline{X}_{\rm Hod}(\log D) /
\overline{S}_{\rm Hod}(\log E), T)\rightarrow
\overline{S}_{\rm Hod}(\log E), T)
$$
is a geometric morphism.  The existence of this extension (which over the open
subset $S_{DR}\subset \overline{S}_{\rm Hod}(\log E), T)$ is just the
Gauss-Manin family $Hom (X_{DR}/S_{DR}, T)$) combines the Griffiths
transversality of the Hodge filtration and regularity of the Gauss-Manin
connection.
This is discussed in more detail in \cite{SantaCruz} in the case $T=BG$ or
particularly $BGL(r)$---I just wanted to make the point here that the same
thing goes through for any connected very presentable $n$-stack $T$.

The same thing will work in an analytic setting, but in this case we can use
any good connected analytic $n$-stack $T$ as coefficients.

\end{document}